\documentclass[11pt]{article}
\usepackage{amsfonts}
\usepackage{latexsym,amssymb,amsfonts}
\usepackage{mathrsfs}
\usepackage{graphicx}
\usepackage{psfrag}
\usepackage{amsmath, amsbsy}
\usepackage{amsopn, amstext}
\usepackage{amsmath,multicol,amsopn}
\usepackage{latexsym,amssymb}

\def\sqr#1#2{{\vcenter{\vbox{\hrule height.#2pt
              \hbox{\vrule width.#2pt height#1pt \kern#1pt \vrule width.#2pt}
              \hrule height.#2pt}}}}
\def\signed #1{{\unskip\nobreak\hfil\penalty50
              \hskip2em\hbox{}\nobreak\hfil#1
              \parfillskip=0pt \finalhyphendemerits=0 \par}}
\def\endpf{\signed {$\sqr69$}}

\def\dbR{{\mathop{\rm l\negthinspace R}}}

\def\3n{\negthinspace \negthinspace \negthinspace }
\def\2n{\negthinspace \negthinspace }
\def\1n{\negthinspace }

\def\dbE{{\mathbb{E}}}
\def\dbF{{\mathbb{F}}}

\def\dbN{{\mathbb{N}}}

\def\dbP{{\mathbb{P}}}

\def\dbR{{\mathbb{R}}}

\def\ds{\displaystyle}

\def\={\buildrel \triangle \over =}

\def\resp{{\it resp. }}
\def\a{\alpha}
\def\b{\beta}

\def\d{\delta}
\def\e{\varepsilon}

\def\l{\lambda}

 \def\n{\nabla}

\def\t{\times}
\def\f{\varphi}
\def\th{\theta}
\def\om{\omega}

\def\ns{\noalign{\ss} }

\def\pa{\partial}
\def\G{\Gamma}

\def\L{\Lambda}
\def\Si{\Sigma}

\def\Om{\Omega}

\def\cA{{\cal A}}
\def\cB{{\cal B}}

\def\cF{{\cal F}}

\def\cI{{\cal I}}

\def\cL{{\cal L}}

\def\cX{{\cal X}}

\def\mE{{\mathbb{E}}}

\def\no{\noindent}

\def\ms{\medskip}
\def\bs{\bigskip}
\def\q{\quad}
\def\qq{\qquad}
\def\hb{\hbox}

\def\max{\mathop{\rm max}}
\def\min{\mathop{\rm min}}
\def\exp{\mathop{\rm exp}}

\def\pa{\partial}

\def\cd{\cdot}
\def\cds{\cdots}

\def\div{\hbox{\rm div$\,$}}

\def\dist{\hbox{\rm dist$\,$}}

\def\supp{\hbox{\rm supp$\,$}}

\def\|{\Big |}
\def\({\Big (}
\def\){\Big )}
\def\[{\Big[}
\def\]{\Big]}
\def\be{\begin{equation}}
\def\bel{\begin{equation}\label}
\def\ee{\end{equation}}
\def\bt{\begin{theorem}}
\def\bcd{\begin{condition}}
\def\ecd{\end{condition}}
\def\et{\end{theorem}}
\def\bc{\begin{corollary}}
\def\ec{\end{corollary}}
\def\bde{\begin{definition}}
\def\ede{\end{definition}}
\def\bl{\begin{lemma}}
\def\el{\end{lemma}}
\def\bp{\begin{proposition}}
\def\ep{\end{proposition}}
\def\br{\begin{remark}}
\def\er{\end{remark}}
\def\ba{\begin{array}}
\def\ea{\end{array}}
\def\ed{\end{document}}
\def\ns{\noalign{\ms}}
\def\ds{\displaystyle}

\def\square#1{\vbox{\hrule\hbox{\vrule height#1%
     \kern#1\vrule}\hrule}}
\def\rectangle#1#2{\vbox{\hrule\hbox{\vrule height#1%
     \kern#2\vrule}\hrule}}

\font\tenbb=msbm10 \font\sevenbb=msbm7
\font\fivebb=msbm5

\newfam\bbfam
\scriptscriptfont\bbfam=\fivebb
\textfont\bbfam=\tenbb
\scriptfont\bbfam=\sevenbb

\newtheorem{lemma}{Lemma}[section]
\newtheorem{remark}{Remark}[section]

\newtheorem{theorem}{Theorem}[section]
\newtheorem{corollary}{Corollary}[section]

\newtheorem{definition}{Definition}[section]
\newtheorem{proposition}{Proposition}[section]
\newtheorem{condition}{Condition}[section]

\makeatletter
   
   \@addtoreset{equation}{section}
\makeatother

\allowdisplaybreaks

\begin{document}
\title{\bf Exact Controllability for a Refined Stochastic Wave Equation \ms}

\author{Qi L\"u\thanks{School of Mathematics,
Sichuan University, Chengdu 610064, Sichuan
Province, China. The research of this author is
partially supported by NSF of China under grant
11471231 and the NSFC-CNRS Joint Research
Project under grant 11711530142. {\small\it
E-mail:} {\small\tt lu@scu.edu.cn}.} \q and \q
Xu Zhang\thanks{School of Mathematics, Sichuan
University, Chengdu 610064, Sichuan Province,
China. The research of this author is
partially supported by NSF of China under grant
11221101,  the NSFC-CNRS Joint
Research Project under grant 11711530142, the
PCSIRT under grant IRT$\_$16R53 from the Chinese
Education Ministry. {\small\it E-mail:}
{\small\tt zhang$\_$xu@scu.edu.cn}.}}

\date{}

\maketitle

\begin{abstract}

A widely used stochastic wave equation is the
classical wave equation perturbed by a term of
It\^o's integral. We show that this equation is not exactly controllable even
if the controls are effective everywhere in both
the drift and the diffusion terms and also on
the boundary. In some sense this means that some
key feature has been ignored in this model. Then, based
on a stochastic Newton's law, we propose a
refined stochastic wave equation. By means of a
new global Carleman estimate, we establish the
exact controllability of our stochastic wave
equation with three controls. Moreover, we give
a result about the lack of exact
controllability, which shows that the action of
three controls is necessary. Our analysis
indicates that, at least from the point of view
of control theory, the new stochastic wave
equation introduced in this paper is a more
reasonable model than that in the existing
literatures.

\end{abstract}

\bs

\no{\bf Key Words}.  Exact controllability,
  stochastic wave
equation,   stochastic Newton's law, global Carleman estimate, observability estimate.

\ms

\no{\bf AMS subject classifications}. Primary, 93B05; Secondary, 60H15, 93B07, 35B45.


\section{Introduction}


Let $T > 0$, $G \subset \mathbb{R}^{n}$ ($n \in
\mathbb{N}$) be a given bounded domain with a
$C^{2}$ boundary $\G$. Let  $\G_0$ be a suitably
chosen nonempty subset (to be given later) of
$\G$, and $G_0\subset G$ be a nonempty open
subset. Write
$ Q = (0,T) \t G$, $\Si = (0,T) \t \G$, $\Si_0 =
(0,T) \t \G_0$ and $Q_0=(0,T) \t G_0$. Let
$(\Om,\cF,\dbF,\dbP)$ be a complete filtered
probability space with $\dbF=\{\cF_t\}_{t\ge0}$,
the natural filtration generated by a
one-dimensional  standard Brownian motion
$\{W(t)\}_{t\ge0}$. More notations and
assumptions used below  will be given in Section
\ref{160905s2}.

Consider the following controlled stochastic
wave equation:
\begin{equation}\label{system1}
\left\{
\begin{array}{ll}
\ds dy_t-
\sum_{j,k=1}^n(a^{jk}y_{x_j})_{x_k}dt=(a_1\cd
\nabla y + a_2 y+g_1)dt + (a_3y+g_2)
dW(t)&\mbox{ in }Q,\\
\ns\ds y=\chi_{\Si_0}h&\mbox{ on }\Si,\\
\ns\ds y(0)=y_0,\q y_t(0)=y_1&\mbox{ in }G,
\end{array}
\right.
\end{equation}
where the initial datum $(y_0,y_1)\in
L^2(G)\times H^{-1}(G)$, $y$ is the state, and
$g_1,g_2\in L^\infty_\dbF(0,T;H^{-1}(G))$ and
$h\in L^2_\dbF(0,T;L^2(\G_0))$ are three controls.

The equation \eqref{system1} is introduced to
describe the vibration of strings and membranes
perturbed by random forces, as well as the
propagation of waves in random environment (e.g.
\cite[Chapter 2]{DKMNX}). Let us recall below
the derivation of \eqref{system1} in
one-dimensional spatial domain.

\vspace{0.1cm}

Consider the motion of a strand of DNA. Compared
with its length, the diameter of a DNA molecule
is so small that it can be viewed as a long
elastic string.
Usually, a DNA molecule  floats in fluid. It is
always hit by fluid molecules, just as a
particle of pollen floating in fluid.

Without loss of generality, we assume the mass
of the string per unit length is $1$. Denote by
$y(t,x)$ the displacement of the strand (in $\dbR^3$) at time
$t\in [0,+\infty)$ and position $x\in [0,L]$  for some $L>0$.
There are mainly four kinds of forces acting on
the string: an elastic force $F_1(t,x)$, a
friction force $F_2(t,x)$ due to viscosity of
the fluid, an impact force $F_3(t,x)$ from the
flow of the fluid, and a random impulse
$F_4(t,x)$ from impacts of the fluid's
molecules. By Newton's second law, we have that
\begin{equation}\label{1.17-eq1}
\frac{d^2y(t,x)}{dt^2} = F_1(t,x) + F_2(t,x) +
F_3(t,x) + F_4(t,x).
\end{equation}
Similar to the derivation of the deterministic
wave equation, the elastic force
$F_1(t,x)=y_{xx}(t,x)$. The friction depends on
the nature of the fluid. For a fixed $x$, by the
classical theory of Statistical Mechanics (e.g.
\cite[Chapter VI]{Stowe}), the random impulse
$F_4(t,x)$ at $(t,x)$ can be approximated by a
Gaussian white noise with a given spatial
correlation matrix $k(\cd,\cd,y)$,  depending on
the fluid. More precisely,  for $x_1,\,x_2\in
[0,L]$ and $0\leq s \leq t<+\infty$,
$$
\begin{array}{ll}\ds
\mE(F_4(t,x_1) F_4(s,x_2)^\top) =
k(x_1,x_2,y(t,x_1),y(t,x_2))\d(t-s).
\end{array}
$$
Here $\d(\cd)$ is the usual Dirac delta
function.  Then, the equation \eqref{1.17-eq1}
can be rewritten as the following stochastic
wave equation:
\begin{equation}\label{1.17-eq2}
dy_t(t,x) = y_{xx}(t,x)dt + F_2(t,x)dt +
F_3(t,x)dt + \hat k(x,y(t,x))dW(t).
\end{equation}
Here $\hat k(x,y(t,x))=k(x,x,y(t,x),y(t,x))$.
When $y$ is small, we may assume that $k$ is
linear in $y$, that is, $\hat k(x,y(t,x))=k_1
(t,x)y(t,x)$ for a suitable $k_1(\cd,\cd)$.

Many biological events are related to the motion
of the DNA molecules. Hence, there is a strong
motivation to control its motion. Clearly, one can introduce two kinds of
controls. One is a force
applied on the boundary, to control the
displacement of the strand at the boundary
point, the other is the force acted in the
internal of the strand, which can be put in both
the drift and the diffusion terms. These lead to
a model like the control system \eqref{system1}.

Motivated by the above mentioned practical
problem, we introduce the following notion of
exact controllability for \eqref{system1}.

\begin{definition}\label{exact def}
The system \eqref{system1} is called exactly
controllable at the time $T$ if for any
$(y_0,y_1)\in L^2(G)\times H^{-1}(G)$ and
$(y_0',y_1')\in L^2_{\cF_T}(\Om;L^2(G))\times
L^2_{\cF_T}(\Om;H^{-1}(G))$, one can find a
triple of controls $(g_1,g_2,h)\in
L^2_{\dbF}(0,T;H^{-1}(G))\times
L^2_{\dbF}(0,T;H^{-1}(G))\times
L^2_\dbF(0,T;L^2(\G_{0}))$ such that the
corresponding solution $y$ to the system
\eqref{system1} satisfies that $(y(T),y_t(T)) =
(y_0',y_1')$.
\end{definition}

Since three controls are introduced in
\eqref{system1}, one may guess that the desired
exact controllability should be trivially
correct. To ``justify" this, let us recall that,
in \cite{TZ} the null controllability of the
following stochastic heat equation
 \begin{equation}\label{160909system1}
\left\{
\begin{array}{ll}
\ds dp-
\sum_{j,k=1}^n(a^{jk}p_{x_j})_{x_k}dt=(a_1\cd
\nabla p + a_2 p+\chi_{G_0}(x)u_1)dt +
(a_3p+u_2)
dW(t)&\mbox{ in }Q,\\
\ns\ds p=0&\mbox{ on }\Si,\\
\ns\ds p(0)=p_0&\mbox{ in }G\end{array} \right.
\end{equation}
was achieved by means of two controls $u_1\in
L^2_{\dbF}(0,T;L^2(G_0))$ and $u_2\in
L^2_{\dbF}(0,T;L^2(G))$, where $p_0\in L^2(G)$
is the initial state. Further, one can easily
prove the exact trajectory controllability for
the equation \eqref{160909system1} with the same
type of controls (Note that, exactly for the
same reason as that in the deterministic
setting, one cannot expect the usual exact
controllability for the stochastic heat
equation). On the other hand, in \cite{Lu4,
Luqi8} the exact controllability of stochastic
Schr\"odinger and transport equations were also
obtained by   a boundary control acted on the
drift term (like $h$ in \eqref{system1}) and   a
distributed control imposed on the diffusion
term.

Surprisingly, as we shall show in Theorem
\ref{th-non-con} (in Section 2) that, the exact controllability
of \eqref{system1} fails for any $T>0$ and
$\G_0\subset\G$, even if the controls $g_1$ and
$g_2$ are acted everywhere on the domain $G$ and
$\G_0=\G$. Note that, such kinds of controls are
the strongest control actions that one can
introduce into \eqref{system1}. Obviously, this
differs significantly from the well-known
controllability property of deterministic wave
equations (See \cite{Lions1, Zhangxu6, Zuazua}
and the rich references therein). Since
\eqref{system1} is a generalization of the
classical wave equation to the stochastic
setting, from the viewpoint of control theory,
we believe that some key feature has been ignored in
the derivation of the equation \eqref{system1}.

Motivated by the above-mentioned negative
controllability result for \eqref{system1}, in
what follows, we shall propose a refined model
to describe the DNA molecule. For this purpose,
we partially employ a dynamical theory of
Brownian motions, developed in \cite{Nelson}, to
describe the motion of a particle perturbed by
random forces. In our opinion, the essence of
the theory in \cite{Nelson} is a stochastic
Newton's law, at least in certain suitable
sense.

 According to
\cite[Chapter 11]{Nelson}, we may suppose that
\begin{equation}\label{8.6-eq1.1}
y(t,x)=\int_0^t v(s,x)ds + \int_0^t
F(s,x,y(s))dW(s).
\end{equation}
Here $v(\cd,\cd)$ is the expected  velocity,
$F(\cd,\cd,\cd)$ is the random perturbation from
the fluid molecule. When $y$ is small, one can
assume that $F(\cd,\cd,\cd)$ is linear in the
third argument, i.e.,
\begin{equation}\label{3.18-eq2}
F(s,x,y(t,x))=b_1 (t,x)y(t,x)
\end{equation}
for a suitable $b_1(\cd,\cd)$.

The acceleration at position $x$ along the
string at time $t$ is $v_{t}(t,x)$.  By Newton's
law, it follows that
\begin{equation}\label{1.17-eq1.1}
v_{t}(t,x) = F_1(t,x) + F_2(t,x) + F_3(t,x) +
F_4(t,x).
\end{equation}
Similar to the derivation of \eqref{1.17-eq2},
we have
\begin{equation}\label{1.17-eq2.1}
dv(t,x) = y_{xx}(t,x)dt + F_2(t,x)dt +
F_3(t,x)dt + k_1 (t,x)y(t,x)dW(t).
\end{equation}
Combining \eqref{8.6-eq1.1}, \eqref{3.18-eq2}
and \eqref{1.17-eq2.1}, we obtain that:
\begin{equation}\label{8.6-eq2.1}
\begin{cases}
dy = vdt + b_1 (t,x)y dW(t) &\mbox{ in
}(0,T)\times
(0,L),\\
\ns\ds dv = y_{xx}dt + F_2dt + F_3dt +
k_1(t,x)y(t,x)dW(t) &\mbox{ in }(0,T)\times
(0,L).
\end{cases}
\end{equation}

Stimulated by \eqref{8.6-eq2.1}, we consider the
following controlled stochastic wave-like
equation:
\begin{equation}\label{system2}
\left\{
\begin{array}{ll}
\ds dy= \hat y dt+(a_4y+f)dW(t) &\mbox{ in }Q,\\
\ns\ds d\hat
y-\sum_{j,k=1}^n(a^{jk}y_{x_j})_{x_k}dt=(a_1 \cd
\nabla y + a_2 y + a_5 g)dt + (a_3y+g)
dW(t)&\mbox{ in }Q,\\
\ns\ds y= \chi_{\Si_0}h &\mbox{ on }\Si,\\
\ns\ds y(0)=y_0,\q \hat y(0)=\hat y_0&\mbox{ in
}G.
\end{array}
\right.
\end{equation}
Here $(y_0,\hat y_0)\in L^2(G)\times H^{-1}(G)$,
$(y,\hat y)$ are the state, and $f \in
L^2_\dbF(0,T;L^2(G))$, $g\in L^2_\dbF(0,T;$
$H^{-1}(G))$ and $h\in L^2_\dbF(0,T;L^2(\G_0))$
are three controls.

\begin{remark}\label{rmk1}
We put controls $f$ and $g$ in the diffusion
terms to get the exact controllability. The
first equation in \eqref{system2} can be
regarded as a family of stochastic differential
equations with a parameter $x\in G$. One can put
a control directly in the diffusion term. On the
other hand, the second equation in
\eqref{system2} is a stochastic partial
differential equation. Usually, if we put a
control in the diffusion term, it may affect the
drift term in one way or another.  Here we
assume that the effect is linear and in the form
of ``$a_5gdt$" as that in the second equation of
\eqref{system2}. One may consider more general
cases, say to add a term like ``$a_6fdt$" into
the first equation of \eqref{system2}. However,
a detailed analysis is beyond the scope of this
paper and will be investigated in future.
\end{remark}

\begin{definition}\label{exact def1}
The system \eqref{system2} is called exactly
controllable at time $T$ if for any $(y_0,\hat
y_0)\in L^2(G)\times H^{-1}(G)$ and $(y_1,\hat
y_1)\in L^2_{\cF_T}(\Om;L^2(G))\times
L^2_{\cF_T}(\Om;H^{-1}(G))$, one can find a
triple of controls $(f,g,h)\in
L_{\dbF}^2(0,T;L^2(G))\times
L_{\dbF}^2(0,T;H^{-1}(G))\times L_{\dbF}^2(0,T;
L^2(\G_{0})) $ such that the corresponding
solution $(y,\hat y)$ to \eqref{system2}
satisfies that $(y(T),\hat y(T)) = (y_1,\hat
y_1)$.
\end{definition}

In this paper,  we shall show that \eqref{system2} is
exactly controllable (See Theorem \ref{exact
th}). Hence, from the viewpoint of
controllability theory, the system
\eqref{system2} is a more reasonable model than
\eqref{system1}. Noting that, we also introduce
three controls into \eqref{system2}, which seems
too many. However, we prove that none of these
three controls can be ignored, and moreover the
two internal controls $f$ and $g$ have to be
effective everywhere in the domain $G$ (See
Theorem \ref{th-non-con1}).

There exist many works on controllability of
deterministic partial differential equations
(PDEs for short). Contributions by D. L. Russell
(\cite{Russell}) and by J.-L. Lions
(\cite{Lions1}) are classical in this field.
Some recent progress can be found in
\cite{Coron, Zhangxu6, Zuazua}. In particular,
one may find many works addressing the exact
controllability problems for deterministic wave
equations (See \cite{AI, Bardos-Lebeau-Rauch1,
Lions1, Russell, Zhangxu2, Zuazua} and the rich
reference therein). However, people know very
little about the controllability problems for
stochastic PDEs. In this respect, we refer to
\cite{barbu1, FL1, GCL, Lu0, Luqi8, TZ} for some
results on the controllability of stochastic
parabolic, complex Ginzburg-Landau,
Kuramoto-Sivashinsky, Schr\"odinger and
transport equations. To the best of  our
knowledge, there exists no nontrivial published
result concerning the exact controllability of
stochastic wave equation.

Compared with the deterministic situation, there
are many new difficulties and phenomena appeared
in the study of controllability problems for
stochastic control systems, even for the systems
governed by stochastic (ordinary) differential
equations (e.g. \cite{LZ3, Peng}). For example, it was shown in \cite{LZ3} that
there exist no Kalman-type rank condition
for the null controllability/approximate
controllability for controlled
stochastic differential
equations. People will meet more
obstacles and substantially extra difficulties in the study of controllability
problems for stochastic PDEs. Some of them are
as follows:

\begin{itemize}
  \item Unlike the deterministic PDEs, the solution of a
stochastic PDE is usually nondifferentiable with
respect to the variable with noise (say, the
time variable considered in this paper);

  \item The usual compactness embedding result does not
remain true for the solution spaces related to
stochastic PDEs;

  \item The diffusion term leads some
  difficulties in establishing observability
  estimate;

  \item The most
essential difficulty is that, compared to their
deterministic counterparts, stochastic
PDEs themselves are much
less-understood.
\end{itemize}

Generally speaking, one can find the following
four main methods for solving the exact
controllability problem of deterministic wave
equations:
\begin{itemize}
  \item The first one is based on the Ingham type
inequality (\cite{AI}). This method works well
for wave equations involved in some special
domains, i.e., intervals and rectangles.
However, it seems that it is very hard to be
applied to equations in general domains.

\item The second one is the classical Rellich-type
multiplier approach (\cite{Lions1}). It is used
to treat wave equations with time independent
coefficients. However, it seems that it does not
work for our problem since the coefficients of
lower order terms  are  time dependent.

\item The third one is the microlocal analysis
approach (\cite{Bardos-Lebeau-Rauch1}).  It is
useful to solve controllability problems for
several kinds of PDEs, such as wave equations,
Schr\"{o}dinger equations and plate equations.
Further, it can give sharp sufficient conditions
for the exact controllability of wave equations.
However, there may be lots of obstacles needed
to be surmounted if one wants to utilize this
approach to study stochastic control problems
(see remarks in Section \ref{sec-com} for more
details).

\item The last one is the global Carleman estimate
(\cite{Zhangxu2}). This approach has the
advantage of being more flexible and allowing to
address variable coefficients. Further, it is
robust with respect to the lower order terms and
can be used to get explicit bounds on the
observability constant/control cost in terms of
the potentials entering in it.
\end{itemize}

In recent years, Carleman estimate was also
employed to study the controllability and
observability problems for some stochastic PDEs
(see \cite{Liu, Lu3, Lu4, TZ, Yuan, Zhangxu3}
and the references therein). \emph{Nevertheless,
as that in the deterministic  setting, generally
speaking, Carleman estimate works well only for
single equation rather than system.}

In this paper, we borrow some idea from the
proof of the observability estimate for
stochastic wave equation (see
\cite{Lu3,Zhangxu3} for example). However, since
(\ref{system2}) is a system (of stochastic equations) rather than a single  stochastic
wave equation,
we cannot simply mimic the method in
\cite{Lu3,Zhangxu3} to solve our problem. To
handle these troubles, we have to derive a completely new
pointwise  identity (see Lemma \ref{hyperbolic1} in Section 6).

The rest of this paper is organized as follows.
In Section \ref{160905s2}, we present the main
results. Section \ref{ssec-well} is devoted to
introducing the adjoint system of systems
\eqref{system1} and \eqref{system2}, and proving
a hidden regularity of solutions to this system.
In Section \ref{sec-well}, we establish the
well-posedness  of the systems \eqref{system1}
and \eqref{system2}. In Section
\ref{sec-reduction}, we transform the exact
controllability problem of \eqref{system2} into
the exact controllability problem of a backward
stochastic wave-like equation. Section
\ref{sec-iden} is addressed to a fundamental
identity for stochastic hyperbolic-like
operators. In Section \ref{sec-ob}, we prove an
observability estimate for a stochastic-wave
like equation. Section \ref{sec-con} is devoted
to proofs of the main results. At last, some
further comments and open problems are given in
Section \ref{sec-com}.


\section{Main results}\label{160905s2}


We begin with some notations.

Denote by $\dbE z$ the (mathematical)
expectation of an integrable random variable
$z:(\Omega,\cF,\dbP)\to \dbR$. Let $H$ be a Banach space. Denote by
$L^{2}_{\dbF}(0,T;H)$ the Banach space
consisting of all $H$-valued and $\dbF$-adapted
processes $X(\cdot)$ such that
$\mathbb{E}(|X(\cdot)|^2_{L^2(0,T;H)}) <
\infty$; by $L^{\infty}_{\dbF}(0,T;H)$ the
Banach space consisting of all $H$-valued and
$\dbF$-adapted, essentially bounded processes;
 and by $C_{\dbF}([0,T];L^r(\Omega;H))$
the Banach space consisting of all $H$-valued
and $\dbF$-adapted processes $X(\cdot)$ such
that $X(\cd):[0,T] \to L^r_{\cF_T}(\Omega;H)$ is
continuous  ($r\in[1,\infty]$). Similarly, one
can define $C^{k}_{\dbF}([0,T];L^{r}(\Om;H))$
for any positive integer $k$. All of these
spaces are endowed with their canonical norms.

In this paper,  for simplicity, we use the
notation $ y_{x_j} \= {\partial y(x)}/{\partial
x_j}$, where $x_j$ is the $j$-th coordinate of a
generic point $x=(x_1,\cdots,
x_n)\in\mathbb{R}^{n}$. In a similar manner, we
use notations $z_{x_j}$, $v_{x_j}$, etc. for the
partial derivatives of $z$ and $v$ with respect
to $x_j$. Also, we denote by $\nu(x) =
(\nu^1(x), \cdots, \nu^n(x))$ the unit outward
normal vector of $\G$ at point $x$. In what
follows, we use $C$ to denote a generic positive
constant depending on $T$, $G$ and $\G_0$
(unless otherwise stated), which may vary from
line to line.

Let $(a^{jk})_{1\leq j,k\leq n} \in
C^3(\overline G;\dbR^{n\times n})$ satisfying
that $a^{jk} = a^{kj}$ $(j,k = 1,2,\cdots, n)$
and for some constant $s_0
> 0$,
\begin{equation}\label{bij}
\sum_{j,k=1}^n a^{jk}\xi^{j}\xi^{k}\geq s_0
|\xi|^2, \,\,\,\,\,\,\,\,\,\, \forall\,
(x,\xi)\= (x,\xi^{1}, \cdots, \xi^{n}) \in G \t
\mathbb{R}^{n}.
\end{equation}
Also we assume that
\begin{equation}\label{coef1}
a_1\in
L^\infty_\dbF(0,T;W^{1,\infty}(G;\dbR^n)),\;
a_2,a_3,a_4 \in
L^\infty_\dbF(0,T;L^\infty(G)),\; a_5\in
L^\infty_\dbF(0,T;W^{1,\infty}_0(G)).
\end{equation}

Let us first give the following negative
controllability result for the system
\eqref{system1}.

\begin{theorem}\label{th-non-con}
The system \eqref{system1} is not exactly
controllable for any $T>0$ and $\G_0\subset\G$.
\end{theorem}

Next, we make the following additional
assumptions on the coefficients $(a^{jk})_{1\leq
j,k\leq n}$:
\begin{condition}\label{condition1}
There exists a positive function $\f(\cdot) \in
C^2(\overline{G})$ satisfying that:

 {\rm (1)} For some constant $\mu_0 >
0$, it holds that
\begin{equation}\label{d1}
\begin{array}{ll}\ds
\sum_{j,k=1}^n\Big\{ \sum_{j',k'=1}^n\Big[
2a^{jk'}(a^{j'k}\f_{x_{j'}})_{x_{k'}} -
a^{jk}_{x_{k'}}a^{j'k'}\f_{x_{j'}} \Big]
\Big\}\xi^{j}\xi^{k} \geq \mu_0
\sum_{j,k=1}^n a^{jk}\xi^{j}\xi^{k}, \\
\ns\ds \hspace{7.5cm} \forall\,
(x,\xi^{1},\cdots,\xi^{n}) \in \overline{G}  \t
\mathbb{R}^n.
\end{array}
\end{equation}
{\rm (2)} There is no critical point of
$\f(\cdot)$ in $\overline{G}$, i.e.,
\be\label{d2}\min_{x\in \overline{G} }|\nabla
\f(x)| > 0. \ee
\end{condition}

The set $\G_0$ is as follows: \vspace{-0.32cm}
\begin{equation}\label{def gamma0}
\G_0 \= \Big\{ x\in \G \;\Big|\;
\sum_{j,k=1}^na^{jk}\f_{x_j}(x)\nu^{k}(x) > 0
\Big\}.
\end{equation}

It is easy to check that if $\f(\cdot) $
satisfies Condition \ref{condition1}, then for
any given constants $\a \geq 1$ and $\b \in
\mathbb{R}$,  $\tilde{\f} = \a\f + \b$ still
satisfies Condition  \ref{condition1} with
$\mu_0$ replaced by $\a\mu_0$. Therefore we may
choose $\f$, $\mu_0$, $c_0>0$, $c_1>0$ and $T$
such that  the following  condition holds:
\begin{condition}\label{condition2}
\begin{equation}\label{con2 eq1}
(1). \q\;\; \frac{1}{4}\sum_{j,k=1}^n
a^{jk}(x)\f_{x_j}(x)\f_{x_k}(x) \geq
R^2_1\=\max_{x\in\overline G}\f(x)\geq R_0^2 \=
\min_{x\in\overline G}\f(x),\q \forall x\in
\overline G.\qq\q
\end{equation}
\hspace{0.4cm} {\rm (2)}. \q $T> T_0\=2 R_1.$\\

  {\rm (3)}. \q $\ds\(\frac{2R_1}{T}\)^2<c_1<\frac{2R_1}{T}$ and $\ds c_1<\min\left\{1,\frac{1}{16|a_5|^4_{L^\infty_\dbF(0,T;L^\infty(G))}}\right\}$.\\

\,\,\,\!   {\rm (4)}. \q $\mu_0 - 4c_1 -c_0 >
\sqrt{R_1}$.
\end{condition}

\begin{remark}\label{rm2}
As we have explained, since $\ds\sum_{j,k=1}^n
a^{jk}\f_{x_j}\f_{x_k} >0$, and one can choose
$\mu_0$ in Condition \ref{condition1} large
enough, Condition \ref{condition2} could be
satisfied obviously. We put it here merely to
emphasize the relationship among $0< c_0< c_1
<1$, $\mu_0$ and $T$. In other words, once
Condition \ref{condition1} is fulfilled,
Condition \ref{condition2} can always be
satisfied.

To be more clear, we give an example for the
choice of $\f$ when $(a^{jk})_{1\leq j,k\leq n}$
is the identity matrix. Let
$x_0\in\dbR^n\setminus\overline{G}$ such that
$|x-x_0|\geq 1$ for all $x\in G$ and
$\ds\a_0=\max_{x\in\overline G}|x-x_0|^2$. Then
for all $\a\geq \max\{\a_0,1\}$, \vspace{-0.2cm}
$$
\a\geq \sqrt{\a}\max_{x\in\overline
G}|x-x_0|.\vspace{-0.1cm}
$$
Let $\f(x)=\a|x-x_0|^2$. Then the left hand side
of \eqref{d1} is reduced to\vspace{-0.1cm}
\begin{equation}\label{d1.1}
\sum_{j=1}^n \f_{x_jx_j}\xi_j^2 = 2 \a|\xi|^2,\q
\forall\, (x,\xi^{1},\cdots,\xi^{n}) \in
\overline{G}  \t \mathbb{R}^n.\vspace{-0.1cm}
\end{equation}
And, \eqref{d1} holds with $\mu_0=2\a$. Further,
it is clear that \eqref{d2} is true
and\vspace{-0.1cm}
$$
\frac{1}{4}\sum_{j,k=1}^n
a^{jk}(x)\f_{x_j}(x)\f_{x_k}(x)
=\frac{1}{4}\sum_{j=1}^n \f_{x_j}(x)^2 =
\a^2|x-x_0|^2\geq \max_{x\in\overline
G}\f(x).\vspace{-0.1cm}
$$
Hence, \eqref{con2 eq1} holds. Next, one can
choose $T$ large enough such that the second and
third inequalities in Condition \ref{condition2}
hold and
$c_1<\min\left\{1,\frac{1}{16|a_5|^4_{L^\infty_\dbF(0,T;L^\infty(G))}}\right\}$.
Let $\a\geq \{\a_0,1,4c_1+c_0\}$. Then $\mu_0 -
4c_1 -c_0 > \sqrt{R_1}$.
\end{remark}
\begin{remark}
To ensure that {\rm (3)} in Condition
\ref{condition2} holds, the larger of
$|a_5|_{L^\infty_\dbF(0,T;L^\infty(G))}$ is
given, the smaller of $c_1$ and the longer of
the time $T$ we should choose. This is
reasonable since $a_5$ stands for the effect of
the control in the diffusion term to the drift
term.   One needs  time to get rid of such
effect.
\end{remark}

We  have the following exact controllability
result for the system \eqref{system2}.

\begin{theorem}\label{exact th}
Let Conditions \ref{condition1} and
\ref{condition2} hold. Then, the system
\eqref{system2} is exactly controllable at time
$T$.
\end{theorem}

As mentioned before, we introduce three controls
in the system \eqref{system2}. At a first
glance, it seems unreasonable, especially for
that the controls in the diffusion term of
\eqref{system2} are acted on the whole domain
$G$. One may ask whether localized controls are
enough or the boundary control can be dropped.
However, the answer is negative.  More
precisely, we have the following negative
result.

\begin{theorem}\label{th-non-con1}
The system \eqref{system2} is not exactly
controllable at any time $T>0$ and
$\G_0\subset\G$ provided that  one of the
following three conditions is satisfied:
\begin{itemize}
  \item[\rm 1)] \ $a_4\in C_{\dbF}([0,T];L^\infty(\Omega))$, $G\setminus \overline{G_0}\neq\emptyset$ and $f$ is supported in
$G_0$;
  \item[\rm 2)] \ $a_3\in C_{\dbF}([0,T];L^\infty(\Omega))$, $G\setminus \overline{G_0}\neq\emptyset$ and $g$ is supported in
$G_0$;
  \item[\rm 3)] \ $h=0$.
\end{itemize}
\end{theorem}
\begin{remark}
Although it is necessary to put controls $f$ and
$g$ on the whole domain, one may suspect that
Theorem \ref{exact th} is trivial. For instance,
one may give a possible ``proof" of Theorem
\ref{exact th} as follows:

Choosing $f=-a_4y$ and $g=-a_3y$, then the
system \eqref{system2} becomes
\begin{equation}\label{system5}
\left\{
\begin{array}{ll}
\ds dy= \hat y dt &\mbox{ in }Q,\\
\ns\ds d\hat
y-\sum_{j,k=1}^n(a^{jk}y_{x_j})_{x_k}dt=(a_1 \cd
\nabla y + a_2 y - a_5a_3 y)dt &\mbox{ in }Q,\\
\ns\ds y= \chi_{\Si_0}h &\mbox{ on }\Si,\\
\ns\ds y(0)=y_0,\q \hat y(0)=\hat y_0&\mbox{ in
}G.
\end{array}
\right.
\end{equation}
This is a wave-like equation with random
coefficients. If one regards the sample point
$\om$ as a parameter, then for every given
$\om\in\Om$, there is a control $u(\cd,\cd,\om)$
such that the solution to \eqref{system5}
fulfills $(y(T,x,\om),\hat
y(T,x,\om))=(y_1(x,\om),\hat y_1(x,\om))$. It is
easy to see that the control constructed in this
way belongs to
$L^2_{\cF_T}(\Om;L^2(0,T;L^2(\G_0)))$. However,
we do not know whether it is adapted to the
filtration $\dbF$ or not. If it is not, then it
means to determine the value of the control at
present, one needs to use information in future,
which is inadmissible in the stochastic context.
\end{remark}
%


\section{Backward
stochastic wave equations}\label{ssec-well}

In order to define solutions to both
\eqref{system1} and \eqref{system2} in a
suitable sense, we need to introduce the
following ``reference" equation:
\begin{equation}\label{bsystem1}
\left\{
\begin{array}{ll}
\ds dz=\hat zdt + (b_5 z+Z)dW(t) &\mbox{ in } Q_\tau,\\
\ns\ds d\hat z - \sum_{j,k=1}^n(a^{jk}z_{x_j})_{x_k} dt = (b_1\cd\nabla z + b_2 z + b_3Z+b_4\widehat Z)dt  +  \widehat Z dW(t)  &\mbox{ in } Q_\tau,\\
\ns\ds z = 0 &\mbox{ on }  \Si_\tau,\\
\ns\ds z(\tau) = z^{\tau},\q \hat z(\tau) = \hat
z^{\tau} &\mbox{ in } G,
\end{array}
\right.\vspace{-0.1cm}
\end{equation}
where $\tau\in (0,T]$, $Q_\tau\=(0,\tau)\times
G$, $\Si_\tau\=(0,\tau)\times \G$, $(z^{\tau},
\hat z^{\tau}) \in L^2_{\cF_\tau}(\Om;H_0^1(G)\t
L^2(G))$ and
$$
b_1\!\in\!
L^\infty_\dbF(0,T;W^{1,\infty}(G;\dbR^n)),\;\;
b_i\!\in\! L^\infty_\dbF(0,T;L^\infty(G)),\;\;
i=2,3,4,\q b_5\!\in\!
L^\infty_\dbF(0,T;W^{1,\infty}_0(G)).
$$

For the convenience of the reader, we first
recall the definition of the solution to
\eqref{bsystem1}.
\begin{definition}\label{def bt sol}
A quadruple of stochastic processes $(z,Z,\hat
z,\widehat Z)\in C_{\dbF}([0,\tau];H^1_0(G))\t
L^2_{\dbF}(0,\tau;$\linebreak $H^1_0(G))\times
C_{\dbF}([0,\tau];L^2(G))\t
L^2_{\dbF}(0,\tau;L^2(G))$ is called a weak
solution of the system \eqref{bsystem1} if for
every $\psi\in C_0^\infty(G)$ and a.e.
$(t,\om)\in [0,\tau]\times\Om$, it holds that
\begin{equation}\label{def id2}
z^\tau(x)-z(t,x)=\int_t^\tau\hat z(s,x)ds +
\int_t^\tau
\big[b_5z(s,x)+Z(s,x)\big]dW(s)\vspace{-0.1cm}
\end{equation}
and\vspace{-0.1cm}
\begin{equation}\label{def id2-1}
\begin{array}{ll}
\ds \int_G  \hat z^{\tau}(x)\psi(x)dx - \int_G
\hat z(t,x)\psi(x)dx  + \int_t^\tau\int_G
\sum_{j,k=1}^n a^{jk}(x)z_{x_j}(s,x)
\psi_{x_k}(x)dxds
\\ \ns\ds = \int_t^\tau \int_G
\big[b_1(s,x)\cd\nabla z(s,x) + b_2(s,x)z(s,x)
+ b_3(s,x)Z(s,x)+ b_4(s,x)\widehat
Z(s,x)\big]\psi(x)dxds\\
\ns\ds\q + \int_t^\tau\int_G \widehat Z(s,x)
\psi(x)dxdW(s).
\end{array}
\end{equation}
\end{definition}

\medskip

Let us recall the following well-posedness
result for \eqref{system2} (e.g.
\cite{Al-Hussein1,  Mahmudov1}).

\begin{lemma}\label{well posed1}
For any $(z^{\tau},\hat z^{\tau})\in
L^2_{\cF_\tau}(\Om;H_0^1(G))\times
L^2_{\cF_\tau}(\Om;L^2(G))$, the system
\eqref{bsystem1} admits a unique solution $
(z,Z,\hat z,\widehat Z)$.
Moreover,
\begin{equation}\label{best1}
\begin{array}{ll}\ds
|z|_{C_\dbF([0,\tau];H_0^1(G))} +
|Z|_{L^2_\dbF(0,\tau;H^1_0(G))}+|\hat
z|_{C_\dbF([0,\tau];L^2(G))} + |\widehat
Z|_{L^2_\dbF(0,\tau;L^2(G))}\\
\ns\ds  \leq Ce^{Cr_1}
\big(|z^{\tau}|_{L^2_{\cF_\tau}(\Om;H_0^1(G))}+|\hat
z^{\tau}|_{L^2_{\cF_\tau}(\Om;L^2(G))}\big),
\end{array}
\end{equation}
where\vspace{-0.2cm}
$$ r_1\=
|b_1|^2_{L^\infty_{\dbF}(0,T;W^{1,\infty}(G;\dbR^n))}
+
\sum_{i=2}^4|b_i|^2_{L^\infty_{\dbF}(0,T;L^{\infty}(G))}+|b_5|^2_{L^\infty_{\dbF}(0,T;W^{1,\infty}_0(G))}.
$$
\end{lemma}

We have the following hidden regularity for
solutions to \eqref{bsystem1}.

\begin{proposition}\label{prop-hid}
Let $(z^{\tau},\hat z^{\tau})\in
L^2_{\cF_\tau}(\Om;H_0^1(G))\times
L^2_{\cF_\tau}(\Om;L^2(G))$. Then the solution
$(z,Z,\hat z,$ $\widehat Z)$ of \eqref{bsystem1}
satisfies $\frac{\pa z}{\pa\nu}\big|_{\G}\in
L^2_{\dbF}(0,\tau;L^2(\G))$. Furthermore,
\begin{equation}\label{hid-eq1}
\Big| \frac{\pa z}{\pa\nu}
\Big|_{L^2_\dbF(0,\tau;L^2(\G))} \leq
Ce^{Cr_1}\big(|z^{\tau}|_{L^2_{\cF_\tau}(\Om;H_0^1(G))}+|\hat
z^{\tau}|_{L^2_{\cF_\tau}(\Om;L^2(G))}\big),
\end{equation}
where the constant $C$ is independent of $\tau$.
\end{proposition}

{\it Proof}\,: For any $h\=(h^1,\cds,h^n)\in
C^1(\dbR_t \t \dbR^n_x ; \dbR^n)$, by It\^o's
formula and the first equation of
\eqref{bsystem1}, we have
 $$
 \ba{ll}\ds
\q d(\hat zh\cdot\nabla z) \3n&\ds=d\hat zh\cdot\nabla z+\hat zh_t\cdot\nabla zdt+\hat z h\cdot \nabla dz+d\hat zh\cdot \nabla dz\\
 \ns
&\ds=d\hat zh\cdot\nabla z+\hat zh_t\cdot\nabla zdt+\hat z h\cdot \nabla \big(\hat zdt + ZdW(t)\big)+d\hat zh\cdot \nabla dz\\
 \ns
&\ds=d\hat zh\cdot\nabla z+\hat zh_t\cdot\nabla
zdt+\frac{1}{2}\Big[\div(\hat z^2h)-(\div h)\hat
z^2\Big]+\hat z h\cdot \nabla ZdW(t)+d\hat
zh\cdot \nabla dz.
 \ea
 $$
Hence, similar to the proofs of \cite[Lemma
3.2]{FYZ} and \cite[Proposition 3.2]{Zhangxu3},
it follows from a direct computation that
\begin{equation}\label{8.6-eq3}
\begin{array}{ll}
\ds -\sum_{k=1}^n\Big[ 2(h\cd\nabla
z)\sum_{j=1}^n a^{jk}z_{x_j} + h^k \Big(\hat z^2
- \sum_{i,j=1}^n
a^{ij}z_{x_i} z_{x_j}\Big) \Big]_{x_k}dt\\
\ns =  \ds 2 \Big[\!-\! d(\hat z h\cd \nabla z)
+\Big(d\hat z \!-\! \sum_{j,k=1}^n\!
(a^{jk}z_{x_j})_{x_k}dt\Big) h \cd \nabla z
\!+\! \hat z h_t\cd \n zdt \!-\!
\sum_{i,j,k=1}^n \!a^{ij}z_{x_i} z_{x_k}
h^k_{x_j}dt\Big] \\
\ns \ds\q  - (\div h) \hat z^2dt +
\sum_{j,k=1}^n z_{x_j} z_{x_k} \div(a^{jk}h)dt +
2d\hat z h\cd\nabla dz +2 \hat z h\cd\nabla
(b_5z+Z)dW(t).
\end{array}
\end{equation}
Since $\G\in C^2$, one can find a vector field
$\xi=(\xi^1,\cdots,\xi^n)\in
C^1(\mathbb{R}^n;\mathbb{R}^n)$ such that
$\xi=\nu$ on $\G$ (See \cite[p. 29]{Lions1}).
Setting $h=\xi$ in \eqref{8.6-eq3}, integrating
it in $Q$, and taking expectation on $\Om$, we
get that
\begin{equation}\label{8.6-eq4}
\begin{array}{ll}
\ds -\mE\int_{\Si_\tau}\sum_{k=1}^n\Big[
2(h\cd\nabla z)\sum_{j=1}^n a^{jk}z_{x_j} + h^k
\Big(\hat z^2 - \sum_{i,j=1}^n
a^{ij}z_{x_i} z_{x_j}\Big) \Big]\nu^kd\G dt\\
\ns =  \ds -2 \mE\int_G \hat z^T h\cd \nabla z^T
dx + 2 \mE\int_G \hat z(0) h\cd \nabla
z(0)dx \\
\ns\ds \q + 2
\int_{Q_\tau}\Big[\big(b_1\cd\nabla z + b_2 z +
b_3Z+b_4\widehat Z \big) h \cd \nabla z + \hat z
h_t\cd \n z - \sum_{j,k,l=1}^n a^{jk}z_{x_j}
z_{x_l}
h^l_{x_k}  \\
\ns \ds\qq\qq\q  - (\div h) \hat z^2 +
\sum_{j,k=1}^n z_{x_k} z_{x_j} \div(a^{jk}h) + 2
\widehat Z h\cd\nabla (b_5z+Z)\Big]dxdt \= \cI.
\end{array}
\end{equation}
Noting that $z=0$ on $(0,\tau)\times \G$, we
have
\begin{equation}\label{8.6-eq5}
\begin{array}{ll}
\ds \mE\int_{\Si_\tau}\sum_{k=1}^n\Big[
2(h\cd\nabla z)\sum_{j=1}^n a^{jk}z_{x_j} + h^k
\Big(\hat z^2 - \sum_{i,j=1}^n
a^{ij}z_{x_i} z_{x_j}\Big)\Big]\nu^kd\G dt\\
\ns\ds =  \mE\int_{\Si_\tau}\Big[
2\(h\cd\nu\frac{\pa z}{\pa\nu}\)\sum_{j,k=1}^n
a^{jk}\nu^k\nu^j\frac{\pa z}{\pa\nu} -
\sum_{i,j,k=1}^n
a^{ij}\nu^k\nu^ih^k\nu^j\Big|\frac{\pa
z}{\pa\nu}\Big|^2\Big]d\G dt \\
\ns\ds = \mE\int_{\Si_\tau} \sum_{j,k=1}^n
a^{jk}\nu^k\nu^j\Big|\frac{\pa z}{\pa\nu}\Big|^2
d\G dt \geq s_0\mE\int_{\Si_\tau} \Big|\frac{\pa
z}{\pa\nu}\Big|^2 d\G dt.
\end{array}
\end{equation}
It follows from Lemma \ref{well posed1} that
$$
|\cI|\leq
Ce^{Cr_1}\big(|z^{\tau}|_{L^2_{\cF_\tau}(\Om;H_0^1(G))}+|\hat
z^{\tau}|_{L^2_{\cF_\tau}(\Om;L^2(G))}\big).
$$
This, together with \eqref{8.6-eq4} and
\eqref{8.6-eq5}, implies that \eqref{hid-eq1}
holds.
\endpf

\begin{remark}
Proposition \ref{prop-hid} shows that, solutions
of \eqref{bsystem1} enjoy a better regularity on
the boundary than the one provided by the
classical trace theorem of Sobolev spaces. Such
kind of result  is called a hidden regularity
(of the solution). There are many studies in
this topic for deterministic PDEs (e.g.
\cite{Lions2}).
\end{remark}


\section{Well-posedness of the systems \eqref{system1} and \eqref{system2}}
\label{sec-well}


In this section, we establish the well-posedness
of systems \eqref{system1} and \eqref{system2} .
Throughout this section, $\Gamma_0$ is any fixed
open subset of $\Gamma$, which is not
necessarily given by \eqref{def gamma0}.

Systems \eqref{system1} and \eqref{system2} are
nonhomogeneous boundary value problems. Like the
deterministic ones (e.g. \cite{Lions1, Lions2}),
their solutions are understood in the sense of
transposition solution.

\begin{definition}
A stochastic process $ y\in
C_{\dbF}([0,T];L^2(\Om;L^2(G)))\cap
C_{\dbF}^1([0,T];L^2(\Om;H^{-1}(G))) $ is a
transposition solution to \eqref{system1} if for
any $\tau\in (0,T]$ and $(z^{\tau},\hat
z^{\tau})\in L^2_{\cF_\tau}(\Om;H_0^1(G))\times
L^2_{\cF_\tau}(\Om;L^2(G))$, we have that
\begin{equation}\label{def id}
\begin{array}{ll}
\ds \mE \langle
y_t(\tau),z^{\tau}\rangle_{H^{-1}(G),H^1_0(G)} -
\langle \hat
y_0,z(0)\rangle_{H^{-1}(G),H^1_0(G)}- \mE\langle
y(\tau),\hat z^{\tau}\rangle_{L^2(G)} + \langle
y_0,\hat z(0)\rangle_{L^2(G)}
\\ \ns\ds = \mE\int_0^\tau
\langle g_1,z\rangle_{H^{-1}(G),H_0^1(G)} dt
 + \mE\int_0^\tau
\langle g_2,Z\rangle_{H^{-1}(G),H_0^1(G)}
dt-\mE\int_0^\tau \int_{\G_0}\frac{\pa
z}{\pa\nu}hd\G ds.
\end{array}
\end{equation}
Here $(z,Z,\hat z,\widehat Z)$ solves
\eqref{bsystem1} with
$$
\begin{array}{ll}\ds
b_1 = -a_1,\q b_2 = -\div a_1+a_2,\q b_3=a_3,\q
b_4=0,\q b_5=0.
\end{array}
$$
\end{definition}
\begin{definition}\label{1-def1}
A pair of stochastic processes
$ (y,\hat y) \in
C_{\dbF}([0,T];L^2(\Om;L^2(G)))\,\times\,
C_{\dbF}([0,T];$ $L^2(\Om;$ $H^{-1}(G))) $
is a transposition solution to \eqref{system2}
if for any $\tau\in (0,T]$ and $(z^{\tau},\hat
z^{\tau})\in L^2_{\cF_\tau}(\Om;H_0^1(G))\times
L^2_{\cF_\tau}(\Om;$ $L^2(G))$, we have that
\begin{equation}\label{def id1}
\begin{array}{ll}
\ds \mE \langle \hat
y(\tau),z^{\tau}\rangle_{H^{-1}(G),H^1_0(G)} -
\langle \hat
y_0,z(0)\rangle_{H^{-1}(G),H^1_0(G)}  -
\mE\langle y(\tau),\hat z^{\tau}
\rangle_{L^2(G)} + \langle y_0,\hat
z(0)\rangle_{L^2(G)}
\\ \ns\ds = -
\mE\int_0^\tau\langle f,\widehat Z
\rangle_{L^2(G)} dt + \mE\int_0^\tau \langle
g,Z\rangle_{H^{-1}(G),H_0^1(G)}
dt-\mE\int_0^\tau \int_{\G_0}\frac{\pa
z}{\pa\nu}hd\G ds .
\end{array}
\end{equation}
Here $(z,Z,\hat z,\widehat Z)$ solves
\eqref{bsystem1} with
$$
\begin{array}{ll}\ds
b_1 = -a_1,\q b_2 = -\div a_1+a_2-a_3a_5,\q
b_3=a_3,\q b_4=-a_4, \q b_5=-a_5.
\end{array}
$$
\end{definition}
\begin{remark}
When $h=0$, both systems \eqref{system1} and
\eqref{system2} are homogeneous boundary value
problems. By the classical theory for stochastic
evolution equations, \eqref{system1} and
\eqref{system2} admit respectively a unique weak
solution (e.g. \cite[Chapter 6]{Prato}) $y\in
C_\dbF([0,T];L^2(\Om;\!L^2(G)))\cap
C_\dbF^1([0,T];\!L^2(\Om;\!H^{-1}(G)))$ and
$(y,\hat y)\!\in\! C_\dbF([0,T];L^2(\Om;\!
L^2(G)))\times
C_\dbF([0,T];L^2(\Om;H^{-1}(G)))$. It follows
from It\^o's formula that these solutions are
respectively transposition solutions to
\eqref{system1} and \eqref{system2}. Then, by
the uniqueness of the transposition solution to
\eqref{system1} (\resp \eqref{system2}), we know
that the transposition solution to
\eqref{system1} (\resp \eqref{system2}) is also
the weak solution to \eqref{system1} (\resp
\eqref{system2}).
\end{remark}

We have the following well-posedness result for
\eqref{system2}.

\begin{proposition}\label{well posed1-1}
For each $(y_0,\hat y_0)\in L^2(G)\times
H^{-1}(G)$, the system \eqref{system2} admits a
unique transposition solution $(y,\hat y)$.
Moreover,
\begin{equation}\label{well posed est1}
\begin{array}{ll}\ds
|(y,\hat y)|_{C_{\dbF}([0,T];L^2(\Om;L^2(G)))\times C_{\dbF}([0,T];L^2(\Om;H^{-1}(G)))}\\
\ns\ds \leq Ce^{Cr_2}\big( |y_0|_{L^2(G)} +
|\hat y_0|_{H^{-1}(G)} +
|f|_{L^2_{\dbF}(0,T;L^2(G))} +
|g|_{L^2_{\dbF}(0,T;H^{-1}(G))}+
|h|_{L^2_{\dbF}(0,T;L^2(\G_{0}))}\big).
\end{array}
\end{equation}
Here
\begin{equation}\label{r2}
r_2 \ =
|a_1|^2_{L_{\dbF}^{\infty}(0,T;W^{1,\infty}(G;\mathbb{R}^{n}))}
+
\sum_{k=2}^4|a_k|^2_{L_{\dbF}^{\infty}(0,T;L^{\infty}(G))}+|a_5|^2_{L_{\dbF}^{\infty}(0,T;W_0^{1,\infty}(G))}.
\end{equation}
\end{proposition}

{\it Proof}\,: {\bf Uniqueness}. Assume that
$(y,\hat y)$ and $(\tilde y, \tilde{\hat y})$
are two transposition solutions of
\eqref{system2}. It follows from Definition
\ref{1-def1} that for any $\tau\in (0,T]$ and
$(z^{\tau},\hat z^{\tau})\in
L^2_{\cF_\tau}(\Om;H_0^1(G))\times
L^2_{\cF_\tau}(\Om; L^2(G))$,
\begin{equation}\label{def id1-1}
\begin{array}{ll}
\ds \mE \langle \hat
y(\tau),z^{\tau}\rangle_{H^{-1}(G),H^1_0(G)}\!\!
-\! \mE\langle y(\tau),\hat z^{\tau}
\rangle_{L^2(G)} = \mE \langle \tilde{\hat
y}(\tau),z^{\tau}\rangle_{H^{-1}(G),H^1_0(G)}\!\!
- \!\mE\langle \tilde y(\tau),\hat z^{\tau}
\rangle_{L^2(G)},
\end{array}
\end{equation}
which implies that
$$
\big(\hat y(\tau), y(\tau)\big)=\big(\tilde{\hat
y}(\tau), \tilde y(\tau)\big),\qq
\forall\,\tau\in (0,T].
$$
Hence, $\big(\hat y, y\big)=\big(\tilde{\hat y},
\tilde y\big)$ in
$C_{\dbF}([0,T];L^2(\Om;L^2(G)))\times
C_{\dbF}([0,T];L^2(\Om;H^{-1}(G)))$.

\vspace{0.2cm}

{\bf Existence}. Since $\chi_{\Sigma_0}h\in
L^2_\dbF(0,T;L^2(\G))$, there exists a sequence
$\{h_m\}_{m=1}^\infty\subset C_\dbF^2([0,T];$
$H^{3/2}(\G))$ with $h_m(0)=0$ for all
$m\in\dbN$ such that
\begin{equation}\label{8.31-eq1}
\lim_{m\to\infty} h_m = \chi_{\Sigma_0}h \q
\mbox{ in }L^2_\dbF(0,T;L^2(\G)).
\end{equation}
For each $m\in\dbN$, we can find a $\tilde h_m
\in C_\dbF^2([0,T];H^2(G))$ such that $\tilde
h_m|_{\G}=h_m$ and $\tilde h_m(0)=0$.

Consider the following equation:
\begin{equation}\label{system2n}
\left\{\!
\begin{array}{ll}
\ds d\tilde y_m=(\tilde {\hat y}_m-\tilde h_{m,t})dt+[a_4(\tilde y_m+\tilde h_{m})+f]dW(t) &\mbox{ in }Q,\\
\ns\ds d\tilde {\hat
y}_m\!-\!\sum_{j,k=1}^n\!(a^{jk}\tilde
y_{m,x_j})_{x_k}dt=(a_1 \cd \nabla \tilde y_m
\!+\! a_2 \tilde y_m \!+\! \zeta_m)dt \!+\!
[a_3(\tilde y_m\!+\!\tilde h_{m})+g]
dW(t)&\mbox{ in }Q, \\
\ns\ds \tilde y_m=\tilde{\hat  y}_m=0&\mbox{ on }\Si,\\
\ns\ds \tilde y_m(0)= y_0,\q \tilde {\hat
y}_m(0)=\hat y_0&\mbox{ in }G,
\end{array}
\right.
\end{equation}
where $\ds\zeta_m = \sum_{j,k=1}^n(a^{jk}\tilde
h_{m,x_j})_{x_k} + a_1 \cd \nabla \tilde h_m +
a_2 \tilde h_m$. By the classical theory of
stochastic evolution equations (e.g.
\cite[Chapter 6]{Prato}), the system
\eqref{system2n} admits a unique mild (also
weak) solution $ (\tilde y_m,\tilde{\hat
y}_m)\in C_\dbF([0,T];L^2(\Om;L^2(G)))\times
C_\dbF([0,T];L^2(\Om;H^{-1}(G))). $

Let $y_m=\tilde y_m+\tilde h_m$ and $\hat
y_m=\tilde{\hat y}_m$.  For any
$m_1,m_2\in\dbN$, by It\^o's formula and
integration by parts, we have that
\begin{equation}\label{8.31-eq5}
\begin{array}{ll}
\ds \mE \langle \hat
y_{m_1}(\tau),z^\tau\rangle_{H^{-1}(G),H^1_0(G)}
- \langle \hat
y_0,z(0)\rangle_{H^{-1}(G),H^1_0(G)}  -
\mE\langle y_{m_1}(\tau),\hat z^\tau
\rangle_{L^2(G)} + \langle y_0,\hat
z(0)\rangle_{L^2(G)}
\\ \ns\ds = -
\mE\int_0^\tau\langle f,\widehat
Z\rangle_{L^2(G)}dt + \mE\int_0^\tau \langle
g,Z\rangle_{H^{-1}(G),H_0^1(G)}
dt-\mE\int_0^\tau\int_{\G}\frac{\pa
z}{\pa\nu}h_{m_1}d\G ds
\end{array}
\end{equation}
and
$$
\begin{array}{ll}
\ds \mE \langle \hat
y_{m_2}(\tau),z^\tau\rangle_{H^{-1}(G),H^1_0(G)}
- \langle \hat
y_0,z(0)\rangle_{H^{-1}(G),H^1_0(G)}  -
\mE\langle y_{m_2}(\tau),\hat z^\tau
\rangle_{L^2(G)} + \langle y_0,\hat
z(0)\rangle_{L^2(G)}
\\ \ns\ds = -
\mE\int_0^\tau \langle f,\widehat
Z\rangle_{L^2(G)}dt + \mE\int_0^\tau \langle
g,Z\rangle_{H^{-1}(G),H_0^1(G)} dt
-\mE\int_0^\tau \int_{\G}\frac{\pa
z}{\pa\nu}h_{m_2}d\G ds.
\end{array}
$$
Consequently,
\begin{eqnarray}\label{8.31-eq2}
\3n\3n\mE \langle \hat y_{m_1}(\tau)\!-\!\hat
y_{m_2}(\tau),z^\tau\rangle_{H^{-1}(G),H^1_0(G)}\!
- \!\mE\langle
y_{m_1}(\tau)\!-\!y_{m_2}(\tau),\hat z^\tau
\rangle_{L^2(G)} \!=\!\!
-\mE\!\int_{\Si_{\tau}}\!\!\frac{\pa
z}{\pa\nu}(h_{m_1}\!\!-\!h_{m_2})d\G ds.
\end{eqnarray}

Let us choose $(z^\tau,\hat z^\tau)\in
L^2_{\cF_\tau}(\Om;H_0^1(G))\times
L^2_{\cF_\tau}(\Om;L^2(G))$  such that
$$
|z^\tau|_{L^2_{\cF_\tau}(\Om;H_0^1(G))}=1,\qq
|\hat z^\tau|_{L^2_{\cF_\tau}(\Om;L^2(G))}=1
$$
and
\begin{equation}\label{8.31-eq3}
\begin{array}{ll}\ds
\mE \langle \hat y_{m_1}(\tau)-\hat
y_{m_2}(\tau),z^\tau\rangle_{H^{-1}(G),H^1_0(G)}
- \mE\langle y_{m_1}(\tau)-y_{m_2}(\tau),\hat
z^\tau
\rangle_{L^2(G)}\\
\ns\ds \geq
\frac{1}{2}\big(|y_{m_1}(\tau)-y_{m_2}(\tau)|_{L^2_{\cF_\tau}(\Om;L^2(G))}
+ |\hat y_{m_1}(\tau)-\hat
y_{m_2}(\tau)|_{L^2_{\cF_\tau}(\Om;H^{-1}(G))}\big).
\end{array}
\end{equation}
It follows from \eqref{8.31-eq2},
\eqref{8.31-eq3}  and Proposition \ref{prop-hid}
that
$$
\begin{array}{ll}
\ds
|y_{m_1}(\tau)-y_{m_2}(\tau)|_{L^2_{\cF_\tau}(\Om;L^2(G))}
+ |\hat y_{m_1}(\tau)-\hat
y_{m_2}(\tau)|_{L^2_{\cF_T}(\Om;H^{-1}(G))}
\\
\ns\ds \leq  2\Big| \mE\int_{\Si_\tau}\frac{\pa
z}{\pa\nu}(h_{m_1}-h_{m_2})d\G ds
\Big|\\
\ns\ds \leq C
|h_{m_1}-h_{m_2}|_{L^2_\dbF(0,T;L^2(\G))}
|(z^\tau,\hat
z^\tau)|_{L^2_{\cF_\tau}(\Omega;H^1_0(G))\times
L^2_{\cF_\tau}(\Omega;L^2(G))}\\
\ns\ds \leq C
|h_{m_1}-h_{m_2}|_{L^2_\dbF(0,T;L^2(\G))},
\end{array}
$$
where the constant $C$ is independent of $\tau$.
Consequently, it holds that
$$
|y_{m_1}-y_{m_2}|_{C_\dbF([0,T];L^2(\Om;L^2(G)))}
+ |\hat y_{m_1}-\hat
y_{m_2}|_{C_\dbF([0,T];L^2(\Om;H^{-1}(G)))} \leq
C |h_{m_1}-h_{m_2}|_{L^2_\dbF(0,T;L^2(\G))}.
$$
This concludes that $\{(y_{m},\hat
y_m)\}_{m=1}^\infty$ is a Cauchy sequence in
$C_\dbF([0,T];L^2(\Om;L^2(G)))\times
C_\dbF([0,T];$ $L^2(\Om;H^{-1}(G)))$. Denote by
$(y,\hat y)$ the limit of $\{(y_m,\hat
y_m)\}_{m=1}^\infty$. Letting $m\to \infty$ in
\eqref{8.31-eq5}, we get that
\begin{equation}\label{8.31-eq6}
\begin{array}{ll}
\ds \mE \langle \hat
y(\tau),z^\tau\rangle_{H^{-1}(G),H^1_0(G)} -
\langle \hat
y_0,z(0)\rangle_{H^{-1}(G),H^1_0(G)}  -
\mE\langle y(\tau),\hat z^\tau \rangle_{L^2(G)}
+ \langle y_0,\hat z(0)\rangle_{L^2(G)}
\\ \ns\ds = -
\mE\int_0^\tau \langle f,\widehat
Z\rangle_{L^2(G)}dt + \mE\int_0^\tau \langle
g,Z\rangle_{H^{-1}(G),H_0^1(G)}
dt-\mE\int_0^\tau \int_{\G_0}\frac{\pa
z}{\pa\nu}h d\G ds .
\end{array}
\end{equation}
Thus, $(y,\hat y)$ is a transposition solution
to \eqref{system2}.

Let us choose $(z^\tau,\hat z^\tau)\in
L^2_{\cF_\tau}(\Om;H_0^1(G))\times
L^2_{\cF_\tau}(\Om;L^2(G))$  such that
$$
|z^T|_{L^2_{\cF_\tau}(\Om;H_0^1(G))}=1,\qq |\hat
z^T|_{L^2_{\cF_\tau}(\Om;L^2(G))}=1
$$
and
\begin{equation}\label{8.31-eq7}
\begin{array}{ll}\ds
\mE \langle \hat
y(\tau),z^\tau\rangle_{H^{-1}(G),H^1_0(G)}   -
\mE\langle y(\tau),\hat z^\tau
\rangle_{L^2(G)}\geq
\frac{1}{2}\big(|y(\tau)|_{L^2_{\cF_\tau}(\Om;L^2(G))}
+ |\hat
y(\tau)|_{L^2_{\cF_\tau}(\Om;H^{-1}(G))}\big).
\end{array}
\end{equation}
Combining \eqref{8.31-eq6}, \eqref{8.31-eq7} and
Proposition \ref{prop-hid}, we obtain that
$$
\begin{array}{ll}
\ds  |y(\tau)|_{L^2_{\cF_\tau}(\Om;L^2(G))} +
|\hat y(\tau)|_{L^2_{\cF_T}(\Om;H^{-1}(G))}
\\
\ns\ds \leq 2\(\big|\langle \hat
y_0,z(0)\rangle_{H^{-1}(G),H^1_0(G)}\big| +
\big|\langle y_0,\hat z(0)\rangle_{L^2(G)}\big|+
\Big|\mE\int_0^\tau \langle f,\widehat
Z\rangle_{L^2(G)}dt\Big|\\
\ns\ds\qq + \Big|\mE\int_0^\tau \langle
g,Z\rangle_{H^{-1}(G),H_0^1(G)} dt\Big|+ \Big|
\mE\int_{\Si_\tau}\frac{\pa z}{\pa\nu}hd\G ds
\Big|\) \\
\ns\ds \leq Ce^{Cr_2}\big( |y_0|_{L^2(G)} +
|\hat y_0|_{H^{-1}(G)} +
|f|_{L^2_{\dbF}(0,T;L^2(G))} +
|g|_{L^2_{\dbF}(0,T;H^{-1}(G))}+
|h|_{L^2_{\dbF}(0,T;L^2(\G_{0}))}\big)\\
\ns\ds \qq\qq\times |(z^\tau,\hat
z^\tau)|_{L^2_{\cF_\tau}(\Omega;H^1_0(G))\times
L^2_{\cF_\tau}(\Omega;L^2(G))}\\
\ns\ds \leq Ce^{Cr_2}\big( |y_0|_{L^2(G)} +
|\hat y_0|_{H^{-1}(G)} +
|f|_{L^2_{\dbF}(0,T;L^2(G))} +
|g|_{L^2_{\dbF}(0,T;H^{-1}(G))}+
|h|_{L^2_{\dbF}(0,T;L^2(\G_{0}))}\big),
\end{array}
$$
where the constant $C$ is independent of $\tau$.
Therefore, we have that
$$
\begin{array}{ll}
\ds |y|_{C_\dbF([0,T];L^2(\Om;L^2(G)))} +
|\hat y|_{C_\dbF([0,T];L^2(\Om;H^{-1}(G)))}
\\
\ns\ds  \leq Ce^{Cr_2}\big( |y_0|_{L^2(G)} +
|\hat y_0|_{H^{-1}(G)} +
|f|_{L^2_{\dbF}(0,T;L^2(G))} +
|g|_{L^2_{\dbF}(0,T;H^{-1}(G))}+
|h|_{L^2_{\dbF}(0,T;L^2(\G_{0}))}\big).
\end{array}
$$
This completes the proof of Proposition \ref{well posed1-1}.
\endpf

\ms

Using the same argument as above, we have the
following well-posedness result for
\eqref{system1} (Hence we omit its proof).

\begin{proposition}\label{well posed}
For each $(y_0,\hat y_0)\in L^2(G)\times
H^{-1}(G)$, the system \eqref{system1} admits a
unique transposition solution $y$. Furthermore,
\begin{equation*}\label{well posed est}
\begin{array}{ll}\ds
 |y|_{C_{\dbF}([0,T];L^2(\Om;L^2(G)))\cap C_{\dbF}^1([0,T];L^2(\Om;H^{-1}(G)))}\\
\ns\ds \leq Ce^{Cr_3}\big( |y_0|_{L^2(G)} +
|y_1|_{H^{-1}(G)}  +
|g_1|_{L^2_{\dbF}(0,T;H^{-1}(G))} +
|g_2|_{L^2_{\dbF}(0,T;H^{-1}(G))}+
|h|_{L^2_{\dbF}(0,T;L^2(\G_{0}))}\big).
\end{array}
\end{equation*}
Here\vspace{-0.3cm}
$$
r_3 \ =
|a_1|^2_{L_{\dbF}^{\infty}(0,T;W^{1,\infty}(G;\mathbb{R}^{n}))}
+
\sum_{k=2}^3|a_2|^2_{L_{\dbF}^{\infty}(0,T;L^{\infty}(G))}.
$$
\end{proposition}


\section{A reduction of the exact controllability problem}
\label{sec-reduction}


Definition \ref{1-def1} is a natural
generalization of the transposition solution
from deterministic wave equations to the
stochastic ones. Accordingly,  one has to
establish observability estimates for
\eqref{bsystem1} to get the exact
controllability of \eqref{system2}. But it is
not so easy. In this section, we give a
reduction of exact controllability problems for
these systems, that is, we show that  these
problems can be transformed to exact
controllability problems for backward stochastic
wave equations.

Consider the following controlled backward
stochastic wave equation:
\begin{equation}\label{system4}
\left\{
\begin{array}{ll}
\ds d\mathbf{y}= \hat{\mathbf{y}} dt+(a_4\mathbf{y}+\mathbf{Y})dW(t) &\mbox{ in }Q,\\
\ns\ds
d\hat{\mathbf{y}}-\sum_{j,k=1}^n(a^{jk}\mathbf{y}_{x_j})_{x_k}dt=(a_1
\cd \nabla \mathbf{y} + a_2 \mathbf{y} + a_5
\widehat{\mathbf{Y}})dt +
(a_3\mathbf{y}+\widehat{\mathbf{Y}})
dW(t)&\mbox{ in }Q,\\
\ns\ds \mathbf{y}= \chi_{\Si_0}\mathbf{h} &\mbox{ on }\Si ,\\
\ns\ds \hat{\mathbf{y}}=0&\mbox{ on }\Si,\\
\ns\ds \mathbf{y}(T)=\mathbf{y}^T,\q
\hat{\mathbf{y}}(T)=\hat{\mathbf{y}}^T&\mbox{ in
}G.
\end{array}
\right.
\end{equation}
Here $(\mathbf{y}^T,\hat{\mathbf{y}}^T)\in
L^2_{\cF_T}(\Om;L^2(G))\times
L^2_{\cF_T}(\Om;H^{-1}(G))$,
$(\mathbf{y},\mathbf{Y},\hat{\mathbf{y}},\widehat{\mathbf{Y}})$
are the state and $\mathbf{h}\in L^2_\dbF(0,T;$
$L^2(\G_0))$ is the control.

To define the solution to \eqref{system4}, we
introduce the following (forward) equation:
\begin{equation}\label{bsystem2}
\left\{\!\!
\begin{array}{ll}
\ds d\mathbf{z}=\hat{\mathbf{z}}dt + (\mathbf{f}-a_5 \mathbf{z})dW(t) &\mbox{ in } Q^\tau,\\
\ns\ds d\hat{\mathbf{z}}\! -\!
\sum_{j,k=1}^n\!(a^{jk}\mathbf{z}_{x_j})_{x_k}
dt\!
=\! [-a_1\!\cd\!\nabla \mathbf{z} \!+\! (\!-\div a_1\!+\!a_2\!-\!a_3a_5) \mathbf{z}\! +\! a_3\mathbf{f}\! -\!a_4\hat{\mathbf{f}}]dt\!+\!  \hat{\mathbf{f}} dW(t)  &\mbox{ in } Q^\tau,\\
\ns\ds \mathbf{z} = 0 &\mbox{ on } \Si^\tau,\\
\ns\ds \mathbf{z}(\tau) = \mathbf{z}_\tau,\q
\hat{\mathbf{z}}(\tau) = \hat{\mathbf{z}}_\tau
&\mbox{ in } G.
\end{array}
\right.
\end{equation}
Here $Q^\tau\=(\tau,T)\times G$,
$\Si^\tau\=(\tau,T)\times \G$,
$(\mathbf{z}_\tau,\hat{\mathbf{z}}_\tau)\in
L^2_{\cF_\tau}(\Om;H_0^1(G))\times
L^2_{\cF_\tau}(\Om;L^2(G))$, $\mathbf{f}\in
L^2_\dbF(0,T;$ $H_0^1(G))$ and
$\hat{\mathbf{f}}\in L^2_\dbF(0,T;L^2(G))$.

Let us recall the following well-posedness
result for \eqref{bsystem2} (e.g. \cite[Chapter
6]{Prato}).

\begin{lemma}\label{well posed1-2}
For any
$(\mathbf{z}_\tau,\hat{\mathbf{z}}_\tau)\in
L^2_{\cF_\tau}(\Om;H_0^1(G))\times
L^2_{\cF_\tau}(\Om;L^2(G))$, $\mathbf{f}\in
L^2_\dbF(0,T;H_0^1(G))$ and $\hat{\mathbf{f}}\in
L^2_\dbF(0,T;L^2(G))$, the system
\eqref{bsystem2} admits a unique weak solution
$(\mathbf{z},\hat{\mathbf{z}})\in
C_{\dbF}([\tau,T];H^1_0(G)) \times
C_{\dbF}([\tau,T];L^2(G))$. Moreover,
\begin{equation}\label{best1-1}
\begin{array}{ll}\ds
|\mathbf{z}|_{C_\dbF([\tau,T];H_0^1(G))}
+|\hat{\mathbf{z}}|_{C_\dbF([\tau,T];L^2(G))}\\
\ns\ds  \leq Ce^{Cr_4}
\big(|\mathbf{z}_\tau|_{L^2_{\cF_\tau}(\Om;H_0^1(G))}+|\hat{\mathbf{z}}_\tau|_{L^2_{\cF_\tau}(\Om;L^2(G))}
+|\mathbf{f}|_{L^2_\dbF(0,T;H_0^1(G))}+|\hat{\mathbf{f}}|_{L^2_\dbF(0,T;L^2(G))}\big),
\end{array}
\end{equation}
where \vspace{-0.3cm}
$$
r_4\=
|a_1|^2_{L^\infty_{\dbF}(0,T;W^{1,\infty}(G;\dbR^n))}
+
\sum_{i=2}^5|a_i|^2_{L^\infty_{\dbF}(0,T;L^{\infty}(G))}+
\sum_{i=3}^5|a_i|^4_{L^\infty_{\dbF}(0,T;L^{\infty}(G))}+|a_5|^2_{L^\infty_{\dbF}(0,T;W_0^{1,\infty}(G))}
$$
and the constant $C$ is independent of $\tau$.
\end{lemma}

Similar to the proof of Proposition
\ref{prop-hid}, we have
\begin{proposition}\label{prop-hid1}
The solution $(\mathbf{z},\hat{\mathbf{z}})$ to
\eqref{bsystem2} satisfies $\frac{\pa
\mathbf{z}}{\pa\nu}|_{\G}\in
L^2_{\dbF}(\tau,T;L^2(\G))$. Furthermore,
\begin{equation}\label{hid-eq1-1}
\begin{array}{ll}\ds
\Big| \frac{\pa \mathbf{z}}{\pa\nu}
\Big|_{L^2_\dbF(\tau,T;L^2(\G))} \leq Ce^{Cr_4}
\big(|\mathbf{z}_\tau|_{L^2_{\cF_\tau}(\Om;H_0^1(G))}+|\hat{\mathbf{z}}_\tau|_{L^2_{\cF_\tau}(\Om;L^2(G))}
+|\mathbf{f}|_{L^2_\dbF(0,T;H_0^1(G))}+|\hat{\mathbf{f}}|_{L^2_\dbF(0,T;L^2(G))}\big),
\end{array}
\end{equation}
where the constant $C$ is independent of $\tau$.
\end{proposition}
\begin{definition}\label{1-def3}
A quadruple of stochastic processes $
(\mathbf{y},\mathbf{Y},\hat{\mathbf{y}},\widehat{\mathbf{Y}})
\! \in \! C_{\dbF}([0,T];\!L^2(\Om;\!L^2(G)))
\times L^2_{\dbF}(0,T;$ $L^2(G)) \times
C_{\dbF}([0,T]; L^2(\Om; H^{-1}(G))) \times
L^2_{\dbF}(0,T;H^{-1}(G))$
is a transposition solution to \eqref{system4}
if for every
$(\mathbf{z}_\tau,\hat{\mathbf{z}}_\tau)\in
L^2_{\cF_\tau}(\Om;H_0^1(G))\times
L^2_{\cF_\tau}(\Om;L^2(G))$, $\mathbf{f}\in
L^2_\dbF(0,T;H_0^1(G))$ and $\hat{\mathbf{f}}\in
L^2_\dbF(0,T;$ $L^2(G))$, one has that
\begin{eqnarray}\label{def id-eq1}
&&\mE \langle
\hat{\mathbf{y}}^T\!,\mathbf{z}(T)\rangle_{H^{-1}(G),H^1_0(G)}\!
- \mE\langle
\hat{\mathbf{y}}(\tau),\mathbf{z}_\tau\rangle_{H^{-1}(G),H^1_0(G)}\!
- \mE\langle \mathbf{y}^T\!,\hat{\mathbf{z}}(T)
\rangle_{L^2(G)}\! + \mE\langle
\mathbf{y}(\tau),
\hat{\mathbf{z}}_\tau\rangle_{L^2(G)} \nonumber
\\ &&  = -
\mE\int_\tau^T\langle
\mathbf{Y},\hat{\mathbf{f}} \rangle_{L^2(G)} dt
+ \mE\int_\tau^T \langle
\widehat{\mathbf{Y}},\mathbf{f}\rangle_{H^{-1}(G),H_0^1(G)}
dt-\mE\int_\tau^T \int_{\G_0}\frac{\pa
\mathbf{z}}{\pa\nu}\mathbf{h}d\G ds .
\end{eqnarray}
Here $(\mathbf{z},\hat{\mathbf{z}})$ solves
\eqref{bsystem2}.
\end{definition}

We have the following result:
\begin{proposition}\label{well posed1-2.1}
For each $(\mathbf{y}^T,\hat{\mathbf{y}}^T)\in
L^2_{\cF_T}(\Om;L^2(G))\times
L^2_{\cF_T}(\Om;H^{-1}(G))$, the system
\eqref{system4} admits a unique transposition
solution
$(\mathbf{y},\mathbf{Y},\hat{\mathbf{y}},\widehat{\mathbf{Y}})$.
Moreover,
\begin{equation}\label{well posed est1.1}
\begin{array}{ll}\ds
|(\mathbf{y},\mathbf{Y},\hat{\mathbf{y}},\widehat{\mathbf{Y}})|_{
C_{\dbF}([0,T];L^2(\Om;L^2(G)))\times
L^2_{\dbF}(0,T;L^2(G))\times
C_{\dbF}([0,T];L^2(\Om;H^{-1}(G)))\times
L^2_{\dbF}(0,T;H^{-1}(G))}\\
\ns\ds \leq Ce^{Cr_2}\big(
|\mathbf{y}^T|_{L^2_{\cF_T}(\Om;L^2(G))} +
|\hat{\mathbf{y}}^T|_{L^2_{\cF_T}(\Om;H^{-1}(G))}
+
|\mathbf{h}|_{L^2_{\dbF}(0,T;L^2(\G_{0}))}\big).
\end{array}
\end{equation}
Here $r_2$ is given by \eqref{r2}.
\end{proposition}

{\it Proof}\,: The proof is similarly to that for Proposition
\ref{well posed1-1}. We give it here for the
convenience of readers.

{\bf Uniqueness}. Assume that
$(\mathbf{y},\mathbf{Y},\hat{\mathbf{y}},\widehat{\mathbf{Y}})$
and
$(\tilde{\mathbf{y}},\widetilde{\mathbf{Y}},\tilde{\hat{\mathbf{y}}},\widetilde{\widehat{\mathbf{Y}}})$
are two transposition solutions of
\eqref{system4}. By Definition
\ref{1-def3}, for any $\tau\in (0,T]$,
$(\mathbf{z}_{\tau},\hat{\mathbf{z}}_{\tau})\in
L^2_{\cF_\tau}(\Om;H_0^1(G))\times
L^2_{\cF_\tau}(\Om; L^2(G))$, $\mathbf{f}\in
L^2_\dbF(0,T;H_0^1(G))$ and $\hat{\mathbf{f}}\in
L^2_\dbF(0,T;L^2(G))$, we have
\begin{equation*}\label{def id1-1.1}
\!\!\begin{array}{ll} \ds \mE \langle
\hat{\mathbf{y}}(\tau),\mathbf{z}_{\tau}\rangle_{H^{-1}(G),H^1_0(G)}
\!- \mE\langle
\mathbf{y}(\tau),\hat{\mathbf{z}}_{\tau}
\rangle_{L^2(G)}\! + \mE\!\int_\tau^T\!\langle
\mathbf{Y},\hat{\mathbf{f}} \rangle_{L^2(G)} dt
- \mE\!\int_\tau^T\! \langle
\widehat{\mathbf{Y}},\mathbf{f}\rangle_{H^{-1}(G),H_0^1(G)}dt
\\ \ns\ds = \mE \langle \tilde{\hat{\mathbf{y}}}(\tau),\mathbf{z}_{\tau}\rangle_{H^{-1}(G),H^1_0(G)}
\!-\! \mE\langle
\tilde{\mathbf{y}}(\tau),\hat{\mathbf{z}}_{\tau}
\rangle_{L^2(G)}\!+\! \mE\!\int_\tau^T\!\langle
\widetilde{\mathbf{Y}},\hat{\mathbf{f}}
\rangle_{L^2(G)} dt\! -\! \mE\!\int_\tau^T\!
\langle
\widetilde{\widehat{\mathbf{Y}}},\mathbf{f}\rangle_{H^{-1}(G),H_0^1(G)}dt,
\end{array}
\end{equation*}
which implies that
$$
\big(\mathbf{y}(\tau),
\hat{\mathbf{y}}(\tau)\big)=\big(\tilde{\mathbf{y}}(\tau),
\tilde{\hat{\mathbf{y}}}(\tau)\big),\qq
\forall\,\tau\in (0,T]
$$
and
$$
(\mathbf{Y},\widehat{\mathbf{Y}})
=(\widetilde{\mathbf{Y}},\widetilde{\widehat{\mathbf{Y}}})
\;\mbox{ in }\; L^2_{\dbF}(0,T;L^2(G)) \times
L^2_{\dbF}(0,T;H^{-1}(G)).
$$
Hence,
$(\mathbf{y},\mathbf{Y},\hat{\mathbf{y}},\widehat{\mathbf{Y}})=(\tilde{\mathbf{y}},\widetilde{\mathbf{Y}},\tilde{\hat{\mathbf{y}}},\widetilde{\widehat{\mathbf{Y}}})$.

\vspace{0.2cm}

{\bf Existence}. Since
$\chi_{\Sigma_0}\mathbf{h}\in
L^2_\dbF(0,T;L^2(\G))$, there exists a sequence
$\{\mathbf{h}_m\}_{m=1}^\infty\subset
C_\dbF^2([0,T];$ $H^{3/2}(\G))$ with
$\mathbf{h}_m(T)=0$ for all $m\in\dbN$  such
that
\begin{equation}\label{8.20-eq1}
\lim_{m\to\infty} \mathbf{h}_m =
\chi_{\Sigma_0}\mathbf{h} \q \mbox{ in
}L^2_\dbF(0,T;L^2(\G)).
\end{equation}
For each $m\in\dbN$, let us choose
$\tilde{\mathbf{h}}_m \in
C_\dbF^2([0,T];H^2(G))$ such that
$\tilde{\mathbf{h}}_m|_{\G}=\mathbf{h}_m$ and
$\tilde{\mathbf{h}}_m(0)=0$.

Consider the following backward stochastic wave
equation:
\begin{equation}\label{system4-1}
\left\{
\begin{array}{ll}
\ds d\tilde{\mathbf{y}}_m=(\tilde{\hat{\mathbf{y}}}_m-\tilde{\mathbf{h}}_{m,t})dt+[a_4(\tilde{\mathbf{y}}_m + \tilde{\mathbf{ h}}_{m})+\widetilde{\mathbf{Y}}_m]dW(t) &\mbox{ in }Q,\\
\ns\ds
d\tilde{\hat{\mathbf{y}}}_m-\sum_{j,k=1}^n(a^{jk}\tilde{\mathbf{y}}_{m,x_j})_{x_k}dt=(a_1
\cd \nabla \tilde{\mathbf{y}}_m + a_2
\tilde{\mathbf{y}}_m + a_5
\widetilde{\widehat{\mathbf{Y}}}_m)dt\\
\ns\ds\qq\qq\qq\qq\qq\qq\;\; +
[a_3(\tilde{\mathbf{ y}}_m+\tilde{
\mathbf{h}}_{m})+\widetilde{\widehat{\mathbf{Y}}}_m
+ \mathbf{\zeta}_m]
dW(t)&\mbox{ in }Q,\\
\ns\ds \tilde{\mathbf{y}}_m= \tilde{\hat{\mathbf{y}}}_m=0  &\mbox{ on }\Si ,\\
\ns\ds \mathbf{y}_m(T)=\mathbf{y}^T,\q
\hat{\mathbf{y}}_m(T)=\hat{\mathbf{y}}^T&\mbox{
in }G.
\end{array}
\right.
\end{equation}
where $\ds\mathbf{\zeta}_m =
\sum_{j,k=1}^n(a^{jk}\tilde{\mathbf{
h}}_{m,x_j})_{x_k} + a_1 \cd \nabla
\tilde{\mathbf{h}}_m + a_2 \tilde{\mathbf{
h}}_m$. By the classical theory of backward
stochastic evolution equations (e.g.
\cite{Al-Hussein1,Mahmudov1}), the system
\eqref{system4-1} admits a unique mild (also
weak) solution
$(\tilde{\mathbf{y}}_m,\widetilde{\mathbf{Y}}_m,\tilde{\hat{\mathbf{y}}}_m,\widetilde{\widehat{\mathbf{Y}}}_m)\in
C_{\dbF}([0,T];L^2(\Om;L^2(G)))\times
L^2_{\dbF}(0,T;L^2(G))\times C_{\dbF}([0,T];$
$L^2(\Om;H^{-1}(G)))\times
L^2_{\dbF}(0,T;H^{-1}(G))$.

Let
$(\mathbf{y}_m,\mathbf{Y}_m,\hat{\mathbf{y}}_m,
\widehat{\mathbf{Y}}_m)=(\tilde{\mathbf{
y}}_m+\tilde{\mathbf{h}}_m,\widetilde{\mathbf{Y}}_m,\tilde{\hat{\mathbf{y}}}_m,\widetilde{\widehat{\mathbf{Y}}}_m)$.
Then
$(\mathbf{y}_m,\mathbf{Y}_m,\hat{\mathbf{y}}_m,
\widehat{\mathbf{Y}}_m)\in
C_{\dbF}([0,T];L^2(\Om;$ $L^2(G)))\times
L^2_{\dbF}(0,T;L^2(G))\times C_{\dbF}([0,T];$
$L^2(\Om;H^{-1}(G)))\times
L^2_{\dbF}(0,T;H^{-1}(G))$. For any
$m_1,m_2\in\dbN$, by It\^o's formula, we have
that
\begin{eqnarray}\label{8.20-eq2}
&&\mE \langle
\hat{\mathbf{y}}^T\!,\mathbf{z}(T)\rangle_{H^{-1}(G),H^1_0(G)}\!\!
-\! \mE\langle
\hat{\mathbf{y}}_{m_1}\!(\tau),\mathbf{z}_\tau\rangle_{H^{-1}(G),H^1_0(G)}\!\!
- \mE\langle \mathbf{y}^T\!,\hat{\mathbf{z}}(T)
\rangle_{L^2(G)}\!   + \!\mE\langle
\mathbf{y}_{m_1}\!(\tau),
\hat{\mathbf{z}}_\tau\rangle_{L^2(G)} \nonumber
\\&& = -
\mE\int_\tau^T\langle
\mathbf{Y}_{m_1},\hat{\mathbf{f}}
\rangle_{L^2(G)} dt + \mE\int_\tau^T \langle
\widehat{\mathbf{Y}}_{m_1},\mathbf{f}\rangle_{H^{-1}(G),H_0^1(G)}
dt-\mE\int_\tau^T \int_{\G_0}\frac{\pa
\mathbf{z}}{\pa\nu}\mathbf{h}_{m_1}d\G ds.
\end{eqnarray}
and
\begin{equation*}\label{8.20-eq3}
\begin{array}{ll} \ds  \mE \langle
\hat{\mathbf{y}}^T\!,\mathbf{z}(T)\rangle_{H^{-1}(G),H^1_0(G)}\!
- \mE\langle
\hat{\mathbf{y}}_{m_2}\!(\tau),\mathbf{z}_\tau\rangle_{H^{-1}(G),H^1_0(G)}\!
- \mE\langle \mathbf{y}^T\!,\hat{\mathbf{z}}(T)
\rangle_{L^2(G)}\! + \mE\langle
\mathbf{y}_{m_2}\!(\tau),
\hat{\mathbf{z}}_\tau\rangle_{L^2(G)}
\\ \ns\ds = -
\mE\int_\tau^T\langle
\mathbf{Y}_{m_2},\hat{\mathbf{f}}
\rangle_{L^2(G)} dt + \mE\int_\tau^T \langle
\widehat{\mathbf{Y}}_{m_2},\mathbf{f}\rangle_{H^{-1}(G),H_0^1(G)}
dt-\mE\int_\tau^T \int_{\G_0}\frac{\pa
\mathbf{z}}{\pa\nu}\mathbf{h}_{m_2}d\G ds .
\end{array}
\end{equation*}
Thus,
\begin{equation}\label{8.20-eq4}
\begin{array}{ll}
\ds \mE \langle \hat{\mathbf{y}}_{m_1}(\tau)-
\hat{\mathbf{y}}_{m_2}(\tau),\mathbf{z}_\tau\rangle_{H^{-1}(G),H^1_0(G)}
- \mE\langle
\mathbf{y}_{m_1}(\tau)-\mathbf{y}_{m_2}(\tau),\hat{\mathbf{z}}_\tau
\rangle_{L^2(G)}\\
\ns\ds\q - \mE\int_\tau^T\langle
\mathbf{Y}_{m_1}-\mathbf{Y}_{m_2},\hat{\mathbf{f}}
\rangle_{L^2(G)} dt + \mE\int_\tau^T \langle
\widehat{\mathbf{Y}}_{m_1}-\widehat{\mathbf{Y}}_{m_2},\mathbf{f}\rangle_{H^{-1}(G),H_0^1(G)}
dt
\\ \ns\ds = -\mE\int_\tau^T \int_{\G_0}\frac{\pa \mathbf{z}}{\pa\nu}(\mathbf{h}_{m_1}-\mathbf{h}_{m_2})d\G ds.
\end{array}
\end{equation}
By \eqref{8.20-eq4}, similar to the proof of
Proposition \ref{well posed1-1}, we can show
that $\{(\mathbf{y}_m,
\hat{\mathbf{y}}_m)\}_{n=1}^\infty$ is a Cauchy
sequence in $C_\dbF([0,T];L^2(\Om;L^2(G)))\times
C_\dbF([0,T];L^2(\Om;H^{-1}(G)))$.

Now we handle $\{(\mathbf{Y}_m,
\widehat{\mathbf{Y}}_m)\}_{n=1}^\infty$. Choose
$\mathbf{f}\in L^2_\dbF(0,T;H_0^1(G))$ and
$\hat{\mathbf{f}}\in L^2_\dbF(0,T;L^2(G))$ such
that
$$
|\mathbf{f}|_{L^2_\dbF(0,T;H_0^1(G))}=1,\qq
|\hat{\mathbf{f}}|_{L^2_\dbF(0,T;L^2(G))}=1
$$
and
\begin{equation}\label{8.20-eq5}
\begin{array}{ll}\ds
-\mE\int_0^T\!\langle
\mathbf{Y}_{m_1}\!\!-\!\mathbf{Y}_{m_2},\hat{\mathbf{f}}
\rangle_{L^2(G)} dt + \mE\int_0^T \!\langle
\widehat{\mathbf{Y}}_{m_1}\!\!-\!\widehat{\mathbf{Y}}_{m_2},\mathbf{f}\rangle_{H^{-1}(G),H_0^1(G)}
dt\\
\ns\ds \geq
\frac{1}{2}\big(|\mathbf{Y}_{m_1}-\mathbf{Y}_{m_2}|_{L^2_{\dbF}(0,T;L^2(G))}
+
|\widehat{\mathbf{Y}}_{m_1}-\widehat{\mathbf{Y}}_{m_2}|_{L^2_{\dbF}(0,T;H^{-1}(G))}\big).
\end{array}
\end{equation}

Let
$(\mathbf{z}_\tau,\hat{\mathbf{z}}_\tau)=(0,0)$
and $\tau=0$  in \eqref{8.20-eq4}. It follows
from \eqref{8.20-eq4}, \eqref{8.20-eq5} and
Proposition \ref{prop-hid1} that
$$
\begin{array}{ll}
\ds
|\mathbf{Y}_{m_1}-\mathbf{Y}_{m_2}|_{L^2_{\dbF}(0,T;L^2(G))}
+
|\widehat{\mathbf{Y}}_{m_1}-\widehat{\mathbf{Y}}_{m_2}|_{L^2_{\dbF}(0,T;H^{-1}(G))}
\\
\ns\ds \leq  2\Big| \mE\int_{\Si_\tau}\frac{\pa
\mathbf{z}}{\pa\nu}(\mathbf{h}_{m_1}-\mathbf{h}_{m_2})d\G
ds
\Big|\\
\ns\ds \leq C
|\mathbf{h}_{m_1}-\mathbf{h}_{m_2}|_{L^2_\dbF(0,T;L^2(\G))}
|(\mathbf{z}_\tau,\hat{\mathbf{z}}_\tau)|_{L^2_{\cF_\tau}(\Omega;H^1_0(G))\times
L^2_{\cF_\tau}(\Omega;L^2(G))}\\
\ns\ds \leq C
|\mathbf{h}_{m_1}-\mathbf{h}_{m_2}|_{L^2_\dbF(0,T;L^2(\G))}.
\end{array}
$$
This implies that
$$
|\mathbf{Y}_{m_1}-\mathbf{Y}_{m_2}|_{L^2_{\dbF}(0,T;L^2(G))}
+
|\widehat{\mathbf{Y}}_{m_1}-\widehat{\mathbf{Y}}_{m_2}|_{L^2_{\dbF}(0,T;H^{-1}(G))}
\leq C
|\mathbf{h}_{m_1}-\mathbf{h}_{m_2}|_{L^2_\dbF(0,T;L^2(\G))}.
$$
Therefore, $\{(\mathbf{Y}_m,
\widehat{\mathbf{Y}}_m)\}_{m=1}^\infty$ is a
Cauchy sequence in $L^2_{\dbF}(0,T;L^2(G))\times
L^2_{\dbF}(0,T;H^{-1}(G))$. Denote by
$(\mathbf{y},\mathbf{Y},\hat{\mathbf{y}},\widehat{\mathbf{Y}})$
the limit of
$\{(\mathbf{y}_m,\mathbf{Y}_m,\hat{\mathbf{y}}_m,
\widehat{\mathbf{Y}}_m)\}_{n=1}^\infty$. By
letting $m_1\to \infty$ in \eqref{8.20-eq2}, we
conclude that
\begin{eqnarray}\label{8.20-eq6}
&&\mE \langle
\hat{\mathbf{y}}^T,\mathbf{z}(T)\rangle_{H^{-1}(G),H^1_0(G)}
\!- \mE\langle
\hat{\mathbf{y}}(\tau),\mathbf{z}_\tau\rangle_{H^{-1}(G),H^1_0(G)}
\!- \mE\langle \mathbf{y}^T,\hat{\mathbf{z}}(T)
\rangle_{L^2(G)}\! + \mE\langle
\mathbf{y}(\tau),
\hat{\mathbf{z}}_\tau\rangle_{L^2(G)}\nonumber
\\ && = -
\mE\int_\tau^T\langle
\mathbf{Y},\hat{\mathbf{f}} \rangle_{L^2(G)} dt
+ \mE\int_\tau^T \langle
\widehat{\mathbf{Y}},\mathbf{f}\rangle_{H^{-1}(G),H_0^1(G)}
dt-\mE\int_\tau^T \int_{\G_0}\frac{\pa
\mathbf{z}}{\pa\nu}\mathbf{h}d\G ds.
\end{eqnarray}
Thus,
$(\mathbf{y},\mathbf{Y},\hat{\mathbf{y}},\widehat{\mathbf{Y}})$
is a solution of \eqref{system2}.

Let us choose
$(\mathbf{z}_\tau,\hat{\mathbf{z}}_\tau)\in
L^2_{\cF_\tau}(\Om;H_0^1(G))\times
L^2_{\cF_\tau}(\Om;L^2(G))$  such that
$$
|\mathbf{z}_\tau|_{L^2_{\cF_\tau}(\Om;H_0^1(G))}=1,\qq
|\hat{\mathbf{z}}_\tau|_{L^2_{\cF_\tau}(\Om;L^2(G))}=1
$$
and
\begin{eqnarray}\label{8.20-eq7}
- \mE\langle
\hat{\mathbf{y}}(\tau),\mathbf{z}_\tau\rangle_{H^{-1}(G),H^1_0(G)}
\!+\! \mE\langle \mathbf{y}(\tau),
\hat{\mathbf{z}}_\tau\rangle_{L^2(G)}\! \geq\!
\frac{1}{2}\big(|\mathbf{y}(\tau)|_{L^2_{\cF_\tau}(\Om;L^2(G))}
\!+\!
|\hat{\mathbf{y}}(\tau)|_{L^2_{\cF_\tau}(\Om;H^{-1}(G))}\big).
\end{eqnarray}
Let $\mathbf{f}=0$ and $\hat{\mathbf{f}}=0$ in
\eqref{8.20-eq6}. Combining \eqref{8.20-eq6},
\eqref{8.20-eq7} and Proposition
\ref{prop-hid1}, we obtain that
$$
\begin{array}{ll}
\ds
|\mathbf{y}(\tau)|_{L^2_{\cF_\tau}(\Om;L^2(G))}
+
|\hat{\mathbf{y}}(\tau)|_{L^2_{\cF_\tau}(\Om;H^{-1}(G))}
\\
\ns\ds \leq 2\(\big|\mE\langle
\hat{\mathbf{y}}^T,\mathbf{z}(T)\rangle_{H^{-1}(G),H^1_0(G)}\big|
+ \big|\mE\langle
\mathbf{y}^T,\hat{\mathbf{z}}(T)
\rangle_{L^2(G)}\big|+ \Big|\mE\int_\tau^T
\int_{\G_0}\frac{\pa
\mathbf{z}}{\pa\nu}\mathbf{h}d\G ds\Big|\) \\
\ns\ds \leq
Ce^{Cr_4}\big(|\mathbf{y}^T|_{L^2_{\cF_T}(\Om;L^2(G))}
+
|\hat{\mathbf{y}}^T|_{L^2_{\cF_T}(\Om;H^{-1}(G))}
+
|\mathbf{h}|_{L^2_{\dbF}(0,T;L^2(\G_{0}))}\big),
\end{array}
$$
where the constant $C$ is independent of $\tau$.
Thus, we find that
$$
\begin{array}{ll}
\ds
|\mathbf{y}|_{C_\dbF([0,T];L^2(\Om;L^2(G)))} +
|\hat{\mathbf{y}}|_{C_\dbF([0,T];L^2(\Om;H^{-1}(G)))}
\\
\ns\ds  \leq
Ce^{Cr_4}\big(|\mathbf{y}^T|_{L^2_{\cF_T}(\Om;L^2(G))}
+
|\hat{\mathbf{y}}^T|_{L^2_{\cF_T}(\Om;H^{-1}(G))}
+
|\mathbf{h}|_{L^2_{\dbF}(0,T;L^2(\G_{0}))}\big).
\end{array}
$$
Let us choose $\mathbf{f}\in
L^2_\dbF(0,T;H_0^1(G))$ and $\hat{\mathbf{f}}\in
L^2_\dbF(0,T;L^2(G))$  such that
$$
|\mathbf{f}|_{L^2_\dbF(0,T;H_0^1(G))}=1,\qq
|\hat{\mathbf{f}}|_{L^2_\dbF(0,T;L^2(G))}=1
$$
and
\begin{eqnarray}\label{8.20-eq8}
\mE\int_\tau^T\!\langle
\mathbf{Y},\hat{\mathbf{f}} \rangle_{L^2(G)} dt
\!-\! \mE\int_\tau^T\! \langle
\widehat{\mathbf{Y}},\mathbf{f}\rangle_{H^{-1}(G),H_0^1(G)}
dt\!\geq\!
\frac{1}{2}\big(|\mathbf{Y}|_{L^2_{\dbF}(0,T;L^2(G))}
\!+\!
|\widehat{\mathbf{Y}}|_{L^2_{\dbF}(0,T;H^{-1}(G))}\big).
\end{eqnarray}
Let
$(\mathbf{z}_\tau,\hat{\mathbf{z}}_\tau)=(0,0)$
and $\tau=0$  in \eqref{8.20-eq6}. It follows
from \eqref{8.20-eq6}, \eqref{8.20-eq8} and
Proposition \ref{prop-hid1} that
$$
\begin{array}{ll}
\ds  |\mathbf{Y}|_{L^2_{\dbF}(0,T;L^2(G))} +
|\widehat{\mathbf{Y}}|_{L^2_{\dbF}(0,T;H^{-1}(G))}
\\
\ns\ds \leq 2\(\big|\mE\langle
\hat{\mathbf{y}}^T,\mathbf{z}(T)\rangle_{H^{-1}(G),H^1_0(G)}\big|
+ \big|\mE\langle
\mathbf{y}^T,\hat{\mathbf{z}}(T)
\rangle_{L^2(G)}\big|+ \Big|\mE\int_\tau^T
\int_{\G_0}\frac{\pa
\mathbf{z}}{\pa\nu}\mathbf{h}d\G ds\Big|\) \\
\ns\ds \leq
Ce^{Cr_4}\big(|\mathbf{y}^T|_{L^2_{\cF_T}(\Om;L^2(G))}
+
|\hat{\mathbf{y}}^T|_{L^2_{\cF_T}(\Om;H^{-1}(G))}
+
|\mathbf{h}|_{L^2_{\dbF}(0,T;L^2(\G_{0}))}\big).
\end{array}
$$
This completes the proof.
\endpf

Now we give the definition of the exact
controllability for the system \eqref{system4}.
\begin{definition}\label{1-def4}
The system \eqref{system4} is called exactly
controllable at time $T$ if for any
$(\mathbf{y}^T,\hat{\mathbf{y}}^T)\in
L^2_{\cF_T}(\Om;L^2(G))\times
L^2_{\cF_T}(\Om;H^{-1}(G))$ and
$(\mathbf{y}_0,\hat{\mathbf{y}}_0)\in
L^2(G)\times H^{-1}(G)$, one can find a control
$\mathbf{h}\in L_{\dbF}^2(0,T; L^2(\G_{0}))$
such that the corresponding solution
$(\mathbf{y},\hat{\mathbf{y}})$ of
\eqref{system4} satisfies that
$(\mathbf{y}(0),\hat{\mathbf{y}}(0)) =
(\mathbf{y}_0,\hat{\mathbf{y}}_0)$.
\end{definition}

It is clear that the following result holds.
\begin{proposition}\label{well posed1-3}
Let $\tau=T$ in \eqref{bsystem1} and $\tau=0$ in
\eqref{bsystem2}. If $(z,Z,\hat z,\widehat Z)$
is a solution of \eqref{bsystem1}, then
$(\mathbf{z},\hat{\mathbf{z}})=(z,\hat z)$ is a
solution of \eqref{bsystem2} with the initial
data
$(\mathbf{z}_0,\hat{\mathbf{z}}_0)=(z(0),\hat
z(0))$ and nonhomogeneous terms
$(\mathbf{f},\hat{\mathbf{f}})=(Z,\widehat Z)$.
On the other hand, if
$(\mathbf{z},\hat{\mathbf{z}})$ is a solution of
\eqref{bsystem2}, then $(z,Z,\hat z,\widehat
Z)=(\mathbf{z},\mathbf{f},\hat{\mathbf{z}},\hat{\mathbf{f}})$
is a solution of \eqref{bsystem1} with the final
data $(z(T),\hat
z(T))=(\mathbf{z}(T),\hat{\mathbf{z}}(T))$.
\end{proposition}

By Proposition \ref{well posed1-3}, and
Definitions \ref{1-def1} and \ref{1-def3}, we
have the following result concerning the
relationship between solutions of
\eqref{system2} and \eqref{system4}.

\begin{proposition}\label{well posed1-4}
If $(y,\hat y)$ is a transposition solution of
\eqref{system2}, then
$(\mathbf{y},\mathbf{Y},\hat{\mathbf{y}},
\widehat{\mathbf{Y}})=(y,f,\hat y,g)$ is a
transposition solution of \eqref{system4} with
the final data
$(\mathbf{y}_T,\hat{\mathbf{y}}_T)=(y(T),\hat
y(T))$. On the other hand, if
$(\mathbf{y},\mathbf{Y},\hat{\mathbf{y}},
\widehat{\mathbf{Y}})$ is a transposition
solution of \eqref{system4}, then $(y,\hat
y)=(\mathbf{y},\hat{\mathbf{y}})$ is a
transposition solution of \eqref{system2} with
the initial data $(y(0),\hat
y(0))=(\mathbf{y}(0),\hat{\mathbf{y}}(0))$ and
the nonhomogeneous terms
$(f,g)=(\mathbf{Y},\widehat{\mathbf{Y}})$.
\end{proposition}

By Propositions \ref{well posed1-3} and
\ref{well posed1-4}, and by borrowing some idea
from \cite{Peng}, we have the following fact:

\begin{proposition}\label{prop-exact}
The system \eqref{system2} is exactly
controllable at time $T$ if and only if the
system \eqref{system4} is exactly controllable
at time $T$.
\end{proposition}
{\it Proof}\,: \textbf{The ``if" part}. Let
$(\mathbf{y}_0,\hat{\mathbf{y}}_0) \in
L^2(G)\times H^{-1}(G)$ and
$(\mathbf{y}^T,\hat{\mathbf{y}}^T) \in
L^2_{\cF_T}(\Om;L^2(G))\times
L^2_{\cF_T}(\Om;H^{-1}(G))$ be arbitrarily
given. Since \eqref{system4} is exactly
controllable at time $T$, there is an
$\mathbf{h}\in L_{\dbF}^2(0,T; L^2(\G_{0}))$
such that the corresponding solution
$(\mathbf{y},\mathbf{Y},\hat{\mathbf{y}},\widehat{\mathbf{Y}})$
of \eqref{system4} satisfies that
$(\mathbf{y}(T),\hat{\mathbf{y}}(T)) =
(\mathbf{y}^T,\hat{\mathbf{y}}^T)$ and
$(\mathbf{y}(0),\hat{\mathbf{y}}(0)) =
(\mathbf{y}_0,\hat{\mathbf{y}}_0)$. Hence,
$(y,\hat y)=(\mathbf{y},\hat{\mathbf{y}})$ is a
solution of \eqref{system2} with a triple of
controls
$(f,g,h)=(\mathbf{Y},\widehat{\mathbf{Y}},\mathbf{h})$
such that $(y(0),\hat
y(0))=(\mathbf{y}_0,\hat{\mathbf{y}}_0)$ and
$(y(T),\hat
y(T))=(\mathbf{y}^T,\hat{\mathbf{y}}^T)$. Hence,
the system \eqref{system2} is exactly
controllable at time $T$.

The proof for \textbf{the ``only if" part} is
similar.
\endpf

\ms

The following result shows that the exact
controllability of \eqref{system4} is equivalent
to an observability estimate of \eqref{bsystem2}
with $\mathbf{f}=\hat{\mathbf{f}}=0$.

\begin{proposition}\label{prop-dual1}
The system \eqref{system4} is exactly
controllable at time $T$ if and only if there is
a constant $C>0$ such that for all
$(\mathbf{z}_0,\hat{\mathbf{z}}_0)\in
H_0^1(G)\times L^2(G)$,  it holds that
\begin{equation}\label{prop-dual1-eq1}
|(\mathbf{z}_0,\hat{\mathbf{z}}_0)|_{H_0^1(G)\times
L^2(G)}^2\leq
C\mE\int_0^T\int_{\G_0}\Big|\frac{\pa
\mathbf{z}}{\pa\nu}\Big|^2d\G dt,
\end{equation}
where $\mathbf{z}$ is the solution of
\eqref{bsystem2} with $\tau = 0$,
$\mathbf{f}=\hat{\mathbf{f}}=0$, $\mathbf{z}(0)
= \mathbf{z}_0$ and
$\hat{\mathbf{z}}(0)=\hat{\mathbf{z}}_0$.
\end{proposition}
\begin{remark}
Compared with \eqref{bsystem1}, \eqref{system4}
is a forward stochastic wave equations.
Generally speaking, it is easier to establish an
observability estimate for \eqref{system4} than
to prove an observability estimate for
\eqref{bsystem1}. This is why we introduce the
reduction in this section.
\end{remark}

{\it Proof of Proposition \ref{prop-dual1}}\,:
We use the classical duality argument, and
divide the proof into two parts.

\ms

\textbf{The ``if" part.} Since the system
\eqref{system4} is linear, we only need to show
that  for any
$(\mathbf{y}_0,\hat{\mathbf{y}}_0)\in
L^2(G)\times H^{-1}(G)$, there is a  control
$ \mathbf{h} \in  L_{\dbF}^2(0,T; L^2(\G_{0})) $
such that the corresponding solution of
\eqref{system4} with
$(\mathbf{y}_T,\hat{\mathbf{y}}_T)=(0,0)$
satisfies that
$(\mathbf{y}(0),\hat{\mathbf{y}}(0)=(\mathbf{y}_0,\hat{\mathbf{y}}_0)$.

Set
$$
\begin{array}{ll}\ds
\cX\=\Big\{-\frac{\pa
\mathbf{z}}{\pa\nu}\Big|_{\G_{0}}\;\Big|\;(\mathbf{z},\hat{\mathbf{z}})\hb{
solves the equation
 }\eqref{bsystem1}\mbox{ with some }
(\mathbf{z}_0,\hat{\mathbf{z}}_0)\in
H_0^1(G)\times L^2(G)\Big\}.
\end{array}
$$
Clearly, $\cX$ is a linear subspace of
$L_{\dbF}^2(0,T;L^2(\G_{0}))$. Define a linear
functional $\cL$ on $\cX$ as follows:
$$
\cL\(-\frac{\pa
\mathbf{z}}{\pa\nu}\Big|_{\G_{0}}\)=-\langle
\hat{\mathbf{y}}_1,\mathbf{z}_0\rangle_{H^{-1}(G),H_0^1(G)}
+ \langle \mathbf{y}_1,
\hat{\mathbf{z}}_0\rangle_{L^2(G)}.
$$
By \eqref{prop-dual1-eq1}, $\cL$ is a bounded
linear functional on $\cX$. By the Hahn-Banach
theorem, it can be extended to be a bounded
linear functional on the space
$L_{\dbF}^2(0,T;L^2(\G_{0})) $. For simplicity,
we still use $\cL$ to denote this extension. By
Riesz's representation theorem, there is an
$\mathbf{h}\in L_{\dbF}^2(0,T; L^2(\G_{0})) $
such that
\begin{equation}\label{con eq0}
-\langle
\hat{\mathbf{y}}_1,\mathbf{z}_0\rangle_{H^{-1}(G),H_0^1(G)}
+ \langle \mathbf{y}_1,
\hat{\mathbf{z}}_0\rangle_{L^2(G)} = -
\mE\int_0^T\int_{\G_{0}} \frac{\pa
\mathbf{z}}{\pa\nu}\mathbf{h} d\G dt .
\end{equation}
We claim that the random field $\mathbf{h}$ is
the desired control. In fact, by the definition
of the solution to \eqref{system4}, we have
\begin{equation}\label{con eq1}
-\langle
\hat{\mathbf{y}}(0),\mathbf{z}_0\rangle_{H^{-1}(G),H_0^1(G)}
+ \langle \mathbf{y}(0),
\hat{\mathbf{z}}_0\rangle_{L^2(G)} =  -
\mE\int_0^T\int_{\G_{0}} \frac{\pa
\mathbf{z}}{\pa\nu}\mathbf{h} d\G dt.
\end{equation}
It follows from   \eqref{con eq0} and \eqref{con
eq1} that
\begin{equation*}\label{con eq3}
-\langle
\hat{\mathbf{y}}_1,\mathbf{z}_0\rangle_{H^{-1}(G),H_0^1(G)}
+ \langle \mathbf{y}_1,
\hat{\mathbf{z}}_0\rangle_{L^2(G)} = -\langle
\hat{\mathbf{y}}(0),\mathbf{z}_0\rangle_{H^{-1}(G),H_0^1(G)}
+ \langle \mathbf{y}(0),
\hat{\mathbf{z}}_0\rangle_{L^2(G)}.
\end{equation*}
Noting that $(\mathbf{z}_0,\hat{\mathbf{z}}_0)$
is an arbitrary element in $H_0^1(G)\times
L^2(G)$, we obtain
$(\mathbf{y}(0),\hat{\mathbf{y}}(0))=(\mathbf{y}_1,\hat{\mathbf{y}}_1)$.

\ms

\textbf{The ``only if" part.} We now prove
\eqref{prop-dual1-eq1} by a contradiction
argument.  Otherwise,  one could find a sequence
$\{(\mathbf{z}_{0,k},\hat{\mathbf{z}}_{0,k})\}_{k=1}^\infty\subset
H_0^1(G)\times L^2(G)$ with
$(\mathbf{z}_{0,k},\hat{\mathbf{z}}_{0,k})\neq
(0,0)$ for all $k\in\dbN$, such that the
corresponding solutions
$(\mathbf{z}_k,\hat{\mathbf{z}}_k)$ of
\eqref{bsystem1} with the initial data
$(\mathbf{z}_{0,k},\hat{\mathbf{z}}_{0,k})$
satisfying that
\begin{equation}\label{9.3-eq2}
\begin{array}{ll}\ds
\Big|\frac{\pa
\mathbf{z}_k}{\pa\nu}\Big|_{L_\dbF^2(0,T;L^2(\G_0))}^2
\leq
\frac{1}{k^2}|(\mathbf{z}_{0,k},\hat{\mathbf{z}}_{0,k})|_{H^{1}_0(G)\times
L^2(G)}^2.
\end{array}
\end{equation}
Write
$$
\l_k =
\frac{\sqrt{k}}{|(\mathbf{z}_{0,k},\hat{\mathbf{z}}_{0,k})|_{H^{1}_0(G)\times
L^2(G)}}, \q \tilde{\mathbf{z}}_{0,k} = \l_k
\mathbf{z}_{0,k}, \q \tilde
{\hat{\mathbf{z}}}_{0,k} = \l_k
\hat{\mathbf{z}}_{0,k}
$$
and denote by $(\tilde{\mathbf{z}}_k,
\tilde{\hat{\mathbf{z}}}_k)$ the  solution of
\eqref{bsystem1} (with
$(\mathbf{z}_0,\hat{\mathbf{z}}_0)$ replaced by
$(\tilde{\mathbf{z}}_{k,0},
\tilde{\hat{\mathbf{z}}}_{k,0})$). Then, it
follows from \eqref{9.3-eq2} that, for each $k
\in \dbN$,
\begin{equation}\label{9.3-eq3}
\Big|\frac{\pa
\tilde{\mathbf{z}}_k}{\pa\nu}\Big|_{L_\dbF^2(0,T;L^2(\G_0))}^2
\leq \frac{1}{k}
\end{equation}
and
\begin{equation}\label{9.3-eq4}
|(\tilde{\mathbf{z}}_{k,0},
\tilde{\hat{\mathbf{z}}}_{k,0})|_{
H_0^1(G)\times L^2(G)}=\sqrt{k}.
\end{equation}

Let us choose
$(\mathbf{y}^T,\hat{\mathbf{y}}^T)=(0,0)$ in
\eqref{system4}. Since the system
\eqref{system4} is exactly controllable, for any
given $(\mathbf{y}_1,\hat{\mathbf{y}}_1)\in
L^2(G)\times H^{-1}(G)$, there is a control
$\mathbf{h}\in L_{\dbF}^2(0,T; L^2(\G_{0}))$
driving the corresponding solution of
\eqref{system4} to
$(\mathbf{y}_1,\hat{\mathbf{y}}_1)$. It follows
from the definition of the solution to
\eqref{system4} that
\begin{equation*}\label{9.3-eq1}
-\langle
\hat{\mathbf{y}}_1,\mathbf{z}_0\rangle_{H^{-1}(G),H_0^1(G)}
+ \langle \mathbf{y}_1,
\hat{\mathbf{z}}_0\rangle_{L^2(G)} =
-\mE\int_0^T \int_{\G_0}\frac{\pa
\mathbf{z}}{\pa\nu}\mathbf{h} d\G ds .
\end{equation*}
Thus, for every $k\in\dbN$, we have that
\begin{equation}\label{9.3-eq5}
-\langle
\hat{\mathbf{y}}_1,\tilde{\mathbf{z}}_{k,0}\rangle_{H^{-1}(G),H_0^1(G)}
+ \langle \mathbf{y}_1,
\tilde{\hat{\mathbf{z}}}_{k,0}\rangle_{L^2(G)} =
- \int_0^T \int_{\G_0}\frac{\pa
\tilde{\mathbf{z}}_k}{\pa\nu}\mathbf{h} d\G ds.
\end{equation}
This, together with \eqref{9.3-eq3} and  the
arbitrariness of
$(\mathbf{y}_1,\hat{\mathbf{y}}_1)$, implies
that
$$
(\tilde{\mathbf{z}}_{k,0},
\tilde{\hat{\mathbf{z}}}_{k,0}) \mbox{ tends to
}  0 \mbox{ weakly in } H_0^1(G)\times L^2(G)
\mbox{ as }k\to+\infty.
$$
Hence, by the Principle of Uniform Boundedness,
$\{(\tilde{\mathbf{z}}_{k,0},
\tilde{\hat{\mathbf{z}}}_{k,0})\}_{k=1}^\infty$
is uniformly bounded in $H_0^1(G)\times L^2(G)$,
which contradicts \eqref{9.3-eq4}.
\endpf


\section{A fundamental identity for stochastic hyperbolic-like operators}
\label{sec-iden}


Throughout this section, we assume that $b^{jk}
\in C^{3}([0,T]\t \mathbb{R}^n)$ satisfies
$b^{jk} = b^{kj}$ for  $j,k = 1,2,\cdots,n$, and
$\ell,\Psi \in C^2((0,T)\t\mathbb{R}^n)$. Write
\begin{eqnarray}\label{CH-5-AB1}
\left\{
\begin{array}{lll}\ds
c^{jk}\= (b^{jk}\ell_t)_t +
\sum_{j',k'=1}^n\big[
2b^{jk'}(b^{j'k}\ell_{x_{j'}})_{x_{k'}} -
(b^{jk}b^{j'k'}\ell_{x_{j'}})_{x_{k'}} \big] +
\Psi b^{jk}
\\
\ds\cA\= (\ell_t^2 - \ell_{tt}) -
\sum_{j,k=1}^n\big[b^{jk}\ell_{x_j}\ell_{x_k}
-(b^{jk}\ell_{x_j})_{x_k}\big]-\Psi, \\
\ns\ds \cB\=\cA\Psi +
(\cA\ell_t)_t-\sum_{j,k=1}^n(\cA
b^{jk}\ell_{x_j})_{x_k} + \frac{1}{2}\[
\Psi_{tt} -
\sum_{j,k=1}^n(b^{jk}\Psi_{x_j})_{x_k} \].
\end{array}
\right.
\end{eqnarray}

We shall derive a fundamental identity for the
stochastic hyperbolic-like operator given in the
following  result.

\begin{lemma}\label{hyperbolic1}
Let  $z$ be an $H^2(\mathbb{R}^n)$-valued
semimartingale and $\hat z$ be an
$L^2(\mathbb{R}^n)$-valued semimartingale, and
\bel{160901} dz=\hat zdt + ZdW(t) \qq\mbox{ in }
(0,T)\times \mathbb{R}^n \ee for some $Z\in
L^2_{\dbF}(0,T;H^1(\mathbb{R}^n))$. Set $ \theta
= e^\ell$, $v=\theta z$ and $\hat v = \theta\hat
z + \ell_t v$. Then, for a.e. $x\in
\mathbb{R}^n$ and $\dbP$-a.s. $\omega \in
\Omega$,
\begin{eqnarray}\label{hyperbolic2}
&&\ds\theta \( - 2\ell_t\hat v +
2\sum_{j,k=1}^nb^{jk}\ell_{x_j} v_{x_k}+ \Psi v
\)
\[ d\hat z - \sum_{j,k=1}^n(b^{jk}z_{x_j})_{x_k}  dt \] \nonumber\\
&&\ds\q
+\sum_{j,k=1}^n\[\sum_{j',k'=1}^n\(2b^{jk}
b^{j'k'}\ell_{x_{j'}}v_{x_j}v_{x_{k'}}
-b^{jk}b^{j'k'}\ell_{x_j}v_{x_{j'}}v_{x_{k'}}\)\nonumber\\
&&\ds\q -2\ell_tb^{jk}v_{x_j}\hat
v+b^{jk}\ell_{x_j} \hat v^2 +\Psi
b^{jk}v_{x_j}v-\frac{\Psi_{x_j}}{2}b^{jk}
v^2-\cA b^{jk}\ell_{x_j}
v^2\]_{x_k} \\
&&\ds\q +d
\[\ell_t\sum_{j,k=1}^n b^{jk}v_{x_j}
v_{x_k}+\ell_t\hat v^2-2\sum_{j,k=1}^nb^{jk}\ell_{x_j} v_{x_k} \hat v - \Psi v\hat v+\(\cA\ell_t+\frac{\Psi_t}{2}\)v^2\]\nonumber \\
&&\ds = \Big\{\[\ell_{tt}+\sum_{j,k=1}^n (b^{jk}
\ell_{x_j})_{x_k}-\Psi\]\hat v^2 +\sum_{j,k=1}^n
c^{jk}v_{x_j} v_{x_k}  - 2 \sum_{j,k=1}^n
\((b^{jk}\ell_{x_k})_t +
b^{jk}\ell_{tx_k}\)v_{x_j}\hat
v \nonumber  \\
&&\ds \q    +\cB v^2 \!+\! \(\! - 2\ell_t\hat
v\! +\! 2\!\sum_{j,k=1}^nb^{jk}\ell_{x_j}
v_{x_k}\!+\! \Psi v \)^2\Big\} dt \!+\!
\ell_t(d\hat v)^2 \!-\!
2\!\sum_{j,k=1}^nb^{jk}\ell_{x_j} dv_{x_k} d\hat
v
\!-\! \Psi dvd\hat v \nonumber\\
&&\ds\q +\ell_{t}\sum_{j,k=1}^n b^{jk}(dv_{x_j})
(dv_{x_k})+\cA\ell_t(dv)^2 -\[\theta\(\! - 2\ell_t\hat v \!+\! 2\!\sum_{j,k=1}^nb^{jk}\ell_{x_j} v_{x_k}\!+ \!\Psi v \)\ell_t Z\!\nonumber\\
&&\ds\q -\( 2\!\sum_{j,k=1}^nb^{jk}(\theta
Z)_{x_k}\ell_{x_j}\hat
v\!-\!\theta\Psi_tvZ\!+\!\theta\Psi\hat v Z\)
 +2\(\sum_{j,k=1}^n b^{jk}v_{x_j} (\th Z)_{x_k}\! +\!
\th \cA v  Z\)\ell_t\]dW(t),\nonumber
\end{eqnarray}
where $(dv)^2$  and $(d\hat v)^2$ denote the
quadratic variation processes of $v$ and $\hat
v$, respectively.
\end{lemma}

{\it Proof}\,: By \eqref{160901}, and recalling
$v=\theta z$ and $\hat v =\th \hat z + \ell_t
v$, we obtain that \bel{0707e1} dv=d(\th z) =
\th_t z dt + \th dz = \ell_t \th z dt + \th \hat
zdt + \th  ZdW(t) = \hat v dt+ \th Z dW(t).
 \ee
Hence,
\begin{equation}\label{2c2t1}
\begin{array}{ll}
\displaystyle d\hat z\3n&\ds=d[\th^{-1}(\hat
v-\ell_t v)]=\th^{-1}[d\hat v - \ell_{tt}vdt -
\ell_t dv -
\ell_t (\hat v-\ell_t v)dt]\\
\ns &\ds=\th^{-1} \[d\hat v -\big(2\ell_t\hat
v+\ell_{tt}v-\ell_t^2v\big)dt-\theta\ell_t
ZdW(t)\].
\end{array}
\end{equation}
Similarly, by $b^{jk} = b^{kj}$ for  $j,k =
1,2,\cdots,n$, we have
\begin{equation}\label{w2.1}
\sum_{j,k=1}^n(b^{jk}z_{x_j})_{x_k}
 =\th^{-1}\sum_{j,k=1}^n\big[(b^{jk}v_{x_j})_{x_k}-2b^{jk}\ell_{x_j}v_{x_k}
 +(b^{jk}\ell_{x_j}\ell_{x_k}-b^{jk}_{x_k}\ell_{x_j}
 -b^{jk}\ell_{x_jx_k})v\big].
\end{equation}
Therefore, from \eqref{2c2t1}--\eqref{w2.1} and
the definition of $\cA$ in \eqref{CH-5-AB1}, we
get
\begin{eqnarray}\label{c1e5}
&&\th \( - 2\ell_t\hat v +
2\sum_{j,k=1}^nb^{jk}\ell_{x_j} v_{x_k}+ \Psi v
\)
\[ d\hat z - \sum_{j,k=1}^n(b^{jk}z_{x_j})_{x_k}  dt \]\nonumber\\
&&= \( - 2\ell_t\hat v +
2\sum_{j,k=1}^nb^{jk}\ell_{x_j} v_{x_k}+ \Psi v
\)
\Big[d\hat v-\sum_{j,k=1}^n (b^{jk}v_{x_j})_{x_k}dt+\cA vdt\nonumber\\
&&\q+\( - 2\ell_t\hat v + 2\sum_{j,k=1}^nb^{jk}\ell_{x_j} v_{x_k}+ \Psi v \)dt-\theta\ell_t ZdW(t)\Big]\nonumber\\
&&= \( - 2\ell_t\hat v +
2\sum_{j,k=1}^nb^{jk}\ell_{x_j} v_{x_k}+ \Psi v
\)d\hat v+\( - 2\ell_t\hat v +
2\sum_{j,k=1}^nb^{jk}\ell_{x_j} v_{x_k}+ \Psi v
\)^2dt
\\
&&\q +\( - 2\ell_t\hat v + 2\sum_{j,k=1}^nb^{jk}\ell_{x_j} v_{x_k}+ \Psi v \)\[-\sum_{j,k=1}^n (b^{jk}v_{x_j})_{x_k}+\cA v\]dt\nonumber\\
&&\q-\theta\( - 2\ell_t\hat v +
2\sum_{j,k=1}^nb^{jk}\ell_{x_j} v_{x_k}+ \Psi v
\)\ell_t ZdW(t).\nonumber
\end{eqnarray}

We now analyze the first and third terms in the
right-hand side of~(\ref{c1e5}).

Using It\^o's formula and noting \eqref{0707e1},
we have\vspace{-0.23cm}
\begin{eqnarray}\label{w1}
&&\( - 2\ell_t\hat v + 2\sum_{j,k=1}^nb^{jk}\ell_{x_j} v_{x_k}+ \Psi v \)d\hat v \nonumber\\
&&\ds =d \(-\ell_t\hat
v^2+2\sum_{j,k=1}^nb^{jk}\ell_{x_j} v_{x_k} \hat
v + \Psi v\hat v\)- 2
\sum_{j,k=1}^nb^{jk}\ell_{x_j}
\hat v dv_{x_k} -\Psi \hat vdv\nonumber\\
&& \ds\q -\[\!-\ell_{tt}\hat v^2 \!+ \!2\!\sum_{j,k=1}^n\!\big(b^{jk}\ell_{x_j}\big)_t v_{x_k} \hat v\!+\! \Psi_tv \hat v\]dt\! +\! \ell_t(d\hat v)^2 \!-\! 2\!\sum_{j,k=1}^nb^{jk}\ell_{x_j} dv_{x_k} d\hat v \!-\! \Psi dvd\hat v\nonumber\\
&&\ds =d \(-\ell_t\hat
v^2+2\sum_{j,k=1}^nb^{jk}\ell_{x_j} v_{x_k} \hat
v + \Psi v\hat v-\frac{\Psi_t}{2}v^2\)-
\sum_{j,k=1}^n \big(b^{jk}\ell_{x_j}
\hat v^2 \big)_{x_k}dt\\
&&\ds\q +\Big\{\[\ell_{tt}\!+ \sum_{j,k=1}^n\!
(b^{jk} \ell_{x_j})_{x_k} -\Psi\]\hat v^2  - 2
\sum_{j,k=1}^n\!
(b^{jk}\ell_{x_k})_tv_{x_j}\hat v  +\!\frac{\Psi_{tt}}{2}v^2\Big\}dt +  \ell_t(d\hat v)^2 \nonumber\\
&&\ds \q  - 2\sum_{j,k=1}^nb^{jk}\ell_{x_j}
dv_{x_k} d\hat v - \Psi dvd\hat v   -\[
2\sum_{j,k=1}^nb^{jk}(\theta
Z)_{x_k}\ell_{x_j}\hat
v-\theta\Psi_tvZ+\theta\Psi\hat v
Z\]dW(t).\nonumber
\end{eqnarray}
Next,\vspace{-0.323cm}
\begin{eqnarray}\label{1c1.6}
&&\hspace{-0.55cm}\!\!-2\ell _t\hat v\(-\sum_{j,k=1}^n (b^{jk}v_{x_j})_{x_k}+\cA v\)dt\nonumber\\
&&\hspace{-0.55cm}\!\!= 2\sum_{j,k=1}^n
(\ell_tb^{jk}v_{x_j}\hat
v)_{x_k}dt-2\sum_{j,k=1}^n
\ell_{tx_k}b^{jk}v_{x_j}\hat v
dt-2\ell_t\sum_{j,k=1}^n b^{jk} v_{x_j}\hat
v_{x_k}dt-2\cA\ell_t v\hat
vdt\\
&&\hspace{-0.55cm}\!\!=\! 2\!\sum_{j,k=1}^n\!
(\ell_tb^{jk}v_{x_j}\hat
v)_{x_k}dt\!-\!2\!\sum_{j,k=1}^n\!\!
\ell_{tx_k}b^{jk}v_{x_j}\hat v
dt\!-\!2\ell_t\!\!\sum_{j,k=1}^n\!\! b^{jk}
v_{x_j}(dv\!-\!\th
ZdW(t))_{x_k}\!\!-\!2\cA\ell_t v(dv\!-\!\th
ZdW(t)) \nonumber\\
&&\hspace{-0.55cm}\!\! = 2\!\sum_{j,k=1}^n\!
(\ell_tb^{jk}v_{x_j}\hat
v)_{x_k}dt\!-2\!\sum_{j,k=1}^n\!
\ell_{tx_k}b^{jk}v_{x_j}\hat v
dt\!-d\(\ell_t\!\sum_{j,k=1}^n\! b^{jk}v_{x_j}
v_{x_k}\!+\!\cA\ell_tv^2\)\! + \!\sum_{j,k=1}^n
\big(\ell_tb^{jk}\big)_tv_{x_j}
v_{x_k}dt\nonumber\\
&&\hspace{-0.55cm}\!\!\q+(\cA\ell_t)_tv^2dt
+\ell_{t}\sum_{j,k=1}^n b^{jk}(dv_{x_j})
(dv_{x_k})+\cA\ell_t(dv)^2 +2\[\sum_{j,k=1}^n
b^{jk}v_{x_j} (\th Z)_{x_k} + \th \cA v
Z\]\ell_tdW(t). \nonumber
\end{eqnarray}
Further, by direct computation, one may check
that
\begin{eqnarray}\label{1c1.3}
&&
2\sum_{j,k=1}^nb^{jk}\ell_{x_j}v_{x_k}\(-\sum_{j,k=1}^n (b^{jk}v_{x_j})_{x_k}+\cA v\)\nonumber\\
&& = - \sum_{j,k=1}^n\[\sum_{j',k'=1}^n\(2b^{jk}
b^{j'k'}\ell_{x_{j'}}v_{x_j}v_{x_{k'}}
-b^{jk}b^{j'k'}\ell_{x_j}v_{x_{j'}}v_{x_{k'}}\)-\cA b^{jk}\ell_{x_j} v^2\]_{x_k} \\
&& \q+\sum_{j,k,j',k'=1}^n
\[2b^{jk'}(b^{j'k}\ell_{x_{j'}})_{x_{k'}}-
(b^{jk}b^{j'k'}\ell_{x_{j'}})_{x_{k'}}\]v_{x_j}v_{x_k}-
\sum_{j,k=1}^n(\cA b^{jk}\ell_{x_j})_{x_k}
v^2\nonumber
\end{eqnarray}
and
\begin{equation}\label{2c2t11}
\begin{array}{ll}\ds
\Psi v\(- \sum_{j,k=1}^n
(b^{jk}v_{x_j})_{x_k}+\cA v\)\\
\ns\ds= - \sum_{j,k=1}^n \(\Psi
b^{jk}v_{x_j}v-\frac{\Psi_{x_j}}{2}b^{jk}
v^2\)_{x_k} +\Psi \sum_{j,k=1}^n
b^{jk}v_{x_j}v_{x_k}
+\[-\frac{1}{2}\sum_{j,k=1}^n
(b^{jk}\Psi_{x_j})_{x_k}+\cA\Psi\]v^2.
\end{array}
\end{equation}

Finally, combining \eqref{c1e5}--\eqref{2c2t11},
we arrive at the desired equality
\eqref{hyperbolic2}.
\endpf



\section{Observability estimate for the equation \eqref{bsystem2}}
\label{sec-ob}

In this section, we establish the following
observability estimate for the equation
\eqref{bsystem2}.
\begin{theorem}\label{th-Ob}
Let Conditions \ref{condition1} and
\ref{condition2} hold, and $\Gamma_0$ be given
by \eqref{def gamma0}. Then, all solutions of
the equation \eqref{bsystem2} with
$\mathbf{f}=0$ and $\hat{\mathbf{f}}=0$ satisfy
that
\begin{equation}\label{exact ob est}
|(\mathbf{z}_0,\hat{\mathbf{z}}_0)|_{H^{1}_0(G)\times
L^2(G)} \leq Ce^{Cr_4}\Big|\frac{\pa
\mathbf{z}}{\pa\nu}\Big|_{L_\dbF^2(0,T;L^2(\G_0))}.
\end{equation}
\end{theorem}

{\it Proof} : We borrow some idea from
\cite{Lu3, LZ2}, and divide the proof into three
steps.

\ms

{\bf Step 1.}  Let us choose
\begin{equation}\label{1-eq1}
\ell(t,x)=\l \[\f(x)  -
c_1\(t-\frac{T}{2}\)^2\].
\end{equation}
Put
\begin{equation}\label{1.20-eq2}
\L_{i}\=\Big\{(t,x)\in Q\Big|
\f(x)-c_1\(t-\frac{T}{2}\)
> \frac{R_0^2}{2(i+2)} \Big\}, \q
\mbox{ for }i=0,1,2.
\end{equation}
Let
\begin{equation}\label{1.20-eq3}
\begin{cases}\ds
T_i\=\frac{T}{2}-\e_i T,\qq
T_i'\=\frac{T}{2}+\e_i T,\\
\ns\ds Q_i\= (T_i,T_i')\times G,
\end{cases}
\q\mbox{ for }i=0,1,
\end{equation}
where $\e_0$ and $\e_1$ are given below.

From Condition \ref{condition2} and
\eqref{1-eq1}, we find that
\begin{equation}\label{1.20-eq5}
\ell(0,x)=\ell(T,x)\leq\l\( R_1 -
\frac{c_1T^2}{4}\)<0,\qq \forall x\in G.
\end{equation}
Hence, there exists $\e_1\in (0,\frac{1}{2})$
such that \vspace{-0.3cm}
\begin{equation}\label{1.20-eq6}
\L_{2}\subset Q_1
\end{equation}
and that
\begin{equation}\label{1.20-eq7}
\ell(t,x)<0,\qq \forall\ (t,x)\in [(0,T_1)\cup
(T_1',T)]\times G.
\end{equation}
Next, since $\{T/2\}\times G\subset  \L_0$, we
know that there is an $\e_0>0$ such
that\vspace{-0.13cm}
\begin{equation}\label{1.2-eq5}
Q_0\=\(\frac{T}{2}-\e_0 T,\frac{T}{2}+\e_0
T\)\times G\subset \L_0.
\end{equation}

\medskip

{\bf Step 2.} Apply Lemma \ref{hyperbolic1} with
$(b^{jk})_{1\leq j,k\leq n}=(a^{jk})_{1\leq
j,k\leq n}$ to the solution of the equation
\eqref{bsystem2} with \vspace{-0.13cm}
\begin{equation}\label{1-eq1.1}
\Psi = \ell_{tt} + \sum_{j,k=1}^n
(a^{jk}\ell_{x_j})_{x_k} - c_0\l,
\end{equation}
and then estimate the resulting terms in
\eqref{hyperbolic2} one by one.

Let us first analyze  the terms which stand for
the ``energy" of the solution. To this end, we
need to compute orders of $\l$ in the
coefficients of $\hat v^2$, $|\nabla v|^2$ and
$v^2$.

Clearly, the term for $\hat v^2$ reads
\begin{equation}\label{coeffvt}
\[\ell_{tt}+\sum_{j,k=1}^n (a^{jk}
\ell_{x_j})_{x_k}-\Psi\]\hat v^2 = c_0 \l \hat
v^2.
\end{equation}
Noting that  $\ell_{tx_k}=\ell_{x_kt}=0$, we get
\begin{equation}\label{bcoeffvtvi}
2\sum_{j,k=1}^n
\left[(a^{jk}\ell_{x_k})_t+a^{jk}\ell_{tx_k}\right]v_{x_j}\hat
v = 0.
\end{equation}
From \eqref{1-eq1.1}, we see that
\begin{eqnarray*}
&&\3n\3n\3n
(a^{jk}\ell_t)_t + \sum_{j',k'=1}^n\big[
2a^{jk'}\big(a^{j'k}\ell_{x_{j'}}\big)_{x_{k'}}
- \big(a^{jk}a^{j'k'}\ell_{x_{j'}}\big)_{x_{k'}}
\big] +
\Psi a^{jk} \\
&&\3n\3n\3n\ds =  a^{jk}\ell_{tt}\! +\!
\sum_{j',k'=1}^n\! \big[
2a^{jk'}\big(a^{j'k}\ell_{x_{j'}}\big)_{x_{k'}}
\!-
\big(a^{jk}a^{j'k'}\ell_{x_{j'}}\big)_{x_{k'}}
\big]\! + a^{jk} \[\ell_{tt}\!  +\!
\sum_{j',k'=1}^\infty
\big(a^{j'k'}\ell_{x_{j'}}\big)_{x_{k'}} \! -
c_0\l\]\\
&&\3n\3n\3n\ds = 2 a^{jk}\ell_{tt} +
\sum_{j',k'=1}^n \big[
2a^{jk'}(a^{j'k}\ell_{x_{j'}})_{x_{k'}} -
a^{jk}_{x_{k'}} a^{j'k'}\ell_{x_{j'}} \big] -
a^{jk}c_0\l \\
&&\3n\3n\3n\ds=2 a^{jk}\ell_{tt} \!+  \l\!
\sum_{j',k'=1}^n\! \big[
2a^{jk'}(a^{j'k}\psi_{x_{j'}})_{x_{k'}}\! -\!
a^{jk}_{x_{k'}} a^{j'k'}\psi_{x_{j'}} \big]\!
+\l\! \sum_{j',k'=1}^n\! 2a^{jk'}
a^{j'k}\psi_{x_{j'}}\psi_{x_{k'}}\!  -
a^{jk}c_0\l.
\end{eqnarray*}
This, together with Condition \ref{condition1},
implies that
\begin{equation}\label{vivj}
\begin{array}{ll}
\ds \sum_{j,k=1}^n\!\Big\{ (a^{jk}\ell_t)_t
+\! \sum_{j',k'=1}^n\big[
2a^{jk'}(a^{j'k}\ell_{x_{j'}})_{x_{k'}} \!-
(a^{jk}a^{j'k'}\ell_{x_{j'}})_{x_{k'}} \big]\! +
\Psi a^{jk} \Big\}v_{x_j} v_{x_k} \\
\ns\ds \geq \l\big(\mu_0 - 4c_1 -  c_0
\big)\sum_{j,k=1}^n a^{jk} v_{x_j} v_{x_k}.
\end{array}
\end{equation}
Now we compute the coefficients of $v^2$. By
\eqref{CH-5-AB1}, it is easy to obtain that
\begin{equation}\label{8.21-eq3}
\begin{array}{ll}\ds
\cA \3n &\ds= \ell_t^2 - \ell_{tt} -
\sum_{j,k=1}^n \big[a^{jk}\ell_{x_j} \ell_{x_k}
-
(a^{jk}\ell_{x_j})_{x_k}\big] - \Psi  \\
\ns&\ds   = \l^2 \Big[ c_1^2(2t - T)^2 -
\sum_{j,k=1}^n a^{jk}\f_{x_j} \f_{x_k}  \Big] +
4c_1 \l  +  c_0 \l.
\end{array}
\end{equation}
By the definition of $\cB$, we see that
\begin{equation}\label{B1}
\begin{array}{ll}\ds
\cB \3n& \ds= \cA\Psi + (\cA \ell_t)_t -
\sum_{j,k=1}^n(\cA a^{jk}\ell_{x_j})_{x_k} +
\frac{1}{2}\sum_{j,k=1}^n\big[\Psi_{tt} -
(a^{jk}\Psi_{x_j})_{x_k}\big]
\\
\ns&\ds = 2\cA \ell_{tt} - \l c_0 \cA
-\sum_{j,k=1}^na^{jk}\ell_j\cA_k + \cA_t\ell_t
-\frac{1}{2}\sum_{j,k=1}^n\sum_{j',k'=1}^n\big[
a^{jk}(a^{j'k'}\ell_{x_{j'}})_{x_{k'}x_j}\big]_{x_k} \\
\ns&\ds  = 2\l^3 \Big[-2c_1^3(2t-T)^2 +
2c_1\sum_{j,k=1}^na^{jk}\f_{x_j}\f_{x_k} \Big]
-\l^3c_0c_1^2(2t-T)^2 + \l^3c_0\sum_{j,k=1}^na^{jk}\f_{x_j}\f_{x_k}\\
\ns&\ds \q +
\l^3\sum_{j,k=1}^n\sum_{j',k'=1}^na^{jk}\f_{x_j}(a^{j'k'}\f_{x_{j'}}\f_{x_{k'}})_{x_k}
-
4\l^3c_1^3(2t-T)^2 + O(\l^2)  \\
\ns&\ds  = (4c_1+c_0)\l^3 \sum_{j,k=1}^n
a^{jk}\f_{x_j}\f_{x_k} +
\l^3\sum_{j,k=1}^n\sum_{j',k'=1}^n a^{jk}\f_{x_j}(a^{j'k'}\f_{x_{j'}}\f_{x_{k'}})_{x_k} \\
\ns&\ds \q - (8c_1^3 + c_0c_1^2)\l^3(2t-T)^2 +
O(\l^2).
\end{array}
\end{equation}
Similar to \cite[(3.8)]{Lu3}, we have
\begin{equation}\label{bijdidj}
\sum_{j,k=1}^n\sum_{j',k'=1}^n
a^{jk}\f_{x_j}(a^{j'k'}\f_{x_{j'}}\f_{x_{k'}})_{x_k}\geq
\mu_0\sum_{j,k=1}^n a^{jk}\f_{x_j}\f_{x_k}.
\end{equation}
From \eqref{B1}, \eqref{bijdidj}  and Condition
\ref{condition2}, we obtain that
\begin{equation}\label{B ine}
\begin{array}{ll}\ds   \cB  \geq \l^3
(4c_1\!+\!c_0)\sum_{j,k=1}^n\!a^{jk}\f_{x_j}\f_{x_k}\!
+ \!\l^3 \mu_0 \sum_{j,k=1}^n\!
a^{jk}\f_{x_j}\f_{x_k} \!-\! (8c_1^3
\!+\! 2c_0c_1^2)\l^3(2t\!-\!T)^2 \!+\! O(\l^2)  \\
\ns\ds \q\geq (\mu_0+4c_1+c_0)\l^3
 \sum_{j,k=1}^na^{jk}\f_{x_j}\f_{x_k} -
8c_1^2(4c_1+c_0)\(t-\frac{T}{2}\)^2\l^3 +
O(\l^2).
\end{array}
\end{equation}

Since the diffusion term in the second equation
of \eqref{bsystem2} is zero, we obtain that
\begin{equation}\label{8.21-eq1}
\mE\ell_t(d\hat v)^2=0, \q
\mE\sum_{j,k=1}^na^{jk}\ell_{x_j} dv_{x_k} d\hat
v=0,\q\mE\Psi dvd\hat v=0.
\end{equation}
From \eqref{1-eq1} and noting that
$\mathbf{z}_{x_k}=(\th^{-1}v)_{x_k} =
\th^{-1}v_{x_k}-\th^{-1}\ell_{x_k} v$, we see
that
\begin{eqnarray}\label{8.21-eq2}
&&\hspace{-1.53cm} \ell_{t}\mE\sum_{j,k=1}^n
a^{jk}(dv_{x_j})
(dv_{x_k}) \nonumber\\
&&\hspace{-1.53cm} =
\mE\[2c_1\l^3\(t-\frac{T}{2}\)|a_5|^2\th^2\sum_{j,k=1}^n
a^{jk}\f_{x_j}\f_{x_k}|\mathbf{z}|^2 + 2c_1\l
\(t-\frac{T}{2}\)\th^2|a_5|^2 \sum_{j,k=1}^n
a^{jk}\mathbf{z}_{x_j}\mathbf{z}_{x_k}\nonumber\\
&&\hspace{-1.53cm} \qq + 2c_1\l
\(t-\frac{T}{2}\)\th^2|\mathbf{z}|^2
\sum_{j,k=1}^n a^{jk}a_{5,x_j}a_{5,x_k} +
4c_1\l^2\(t-\frac{T}{2}\)\th^2|a_5|^2 \mathbf{z}
\sum_{j,k=1}^n
a^{jk}\f_{x_j}\mathbf{z}_{x_k}\\
&&\hspace{-1.53cm}\qq +
4c_1\l^2\(t-\frac{T}{2}\)\th^2a_5 |\mathbf{z}|^2
\sum_{j,k=1}^n
a^{jk}\f_{x_j}a_{5,x_k}\] \nonumber\\
&&\hspace{-1.53cm} =
\mE\[4c_1\l^3\(t-\frac{T}{2}\)|a_5|^2
\sum_{j,k=1}^n
a^{jk}\f_{x_j}\f_{x_k}|v|^2+2c_1\l
\(t-\frac{T}{2}\) |a_5|^2 \sum_{j,k=1}^n
a^{jk}v_{x_j}v_{x_k} \nonumber\\
&& \hspace{-1.53cm}\qq   + C(\l^2|\nabla a_5|+
\l |\nabla a_5|^2)|v|^2\].\nonumber
\end{eqnarray}

Next, from \eqref{1-eq1} and \eqref{8.21-eq3},
we find that
\begin{equation}\label{8.21-eq4}
\begin{array}{ll}\ds
\q \mE\big[\cA\ell_t(dv)^2\big]\3n&\ds = \l^3
c_1 (2t - T) \Big[ c_1^2(2t - T)^2 -
\sum_{j,k=1}^n
a^{jk}\f_{x_j} \f_{x_k}  \Big]\mE\big(|a_5|^2 v^2\big)\\
\ns&\ds\q + (4c_1 +  c_0)c_1 (2t - T)
\l^2\mE\big(|a_5|^2 v^2\big).
\end{array}
\end{equation}
From  \eqref{vivj}, \eqref{B ine},
\eqref{8.21-eq2} and \eqref{8.21-eq4}, and
noting the fourth inequality in Condition
\ref{condition2}, we know that there is $c_2>0$
such that for all $(t,x)\in \L_2$, one has that
\begin{eqnarray*}
&&\3n\3n\3n\3n\ds   \mE\Big[\sum_{j,k=1}^n
c^{jk}v_{x_j} v_{x_k} +\cB v^2
+\ell_{t}\sum_{j,k=1}^n a^{jk}(dv_{x_j})
(dv_{x_k})+\cA\ell_t(dv)^2\Big]\nonumber\\
&&\3n\3n\3n\3n\ds =\mE\Big\{\l\[\mu_0 - 4c_1 -
c_0 +2c_1 \(t-\frac{T}{2}\) |a_5|^2
\]\sum_{j,k=1}^n a^{jk} v_{x_j} v_{x_k}+ (\mu_0 + 4c_1 + c_0)\l^3
\sum_{j,k=1}^n\!a^{jk}\f_{x_j}\f_{x_k}|v|^2\\
&& \3n\3n\3n\3n\ds \qq  - 8c_1^2(4c_1 + c_0)\(t
- \frac{T}{2}\)^2\l^3|v|^2  + 4c_1\l^3\(t -
\frac{T}{2}\)|a_5|^2 \sum_{j,k=1}^n
a^{jk}\f_{x_j}\f_{x_k}|v|^2 + O(\l^2)|v|^2\Big\} \nonumber\\
&&\3n\3n\3n\3n\ds \geq \mE\big[c_2\l |\n v|^2 +
c_2\l^3|v|^2+ O(\l^2)|v|^2\big].\nonumber
\end{eqnarray*}
Thus, there exist $\l_1>0$ and $c_3>0$ such that
for all $\l\geq \l_1$ and for every $(t,x)\in
\L_2$, one has that
\begin{equation}\label{8.21-eq5}
\mE\[\sum_{j,k=1}^n c^{jk}v_{x_j} v_{x_k} +\cB
v^2 +\ell_{t}\sum_{j,k=1}^n a^{jk}(dv_{x_j})
(dv_{x_k})+\cA\ell_t(dv)^2\]\geq \mE\big(c_3\l
|\n v|^2 + c_3\l^3|v|^2\big).
\end{equation}
\medskip

{\bf Step 3.} For the boundary terms, by
$v|_{\Sigma}=0$, we have the following equality:
\begin{equation}\label{hyperbolic32}
\begin{array}{ll}\ds
\mathbb{E}\int_{\Si}\sum_{j,k=1}^n\sum_{j',k'=1}^n\Big(
2a^{jk}a^{j'k'}\ell_{x_{j'}}v_{x_j} v_{x_{k'}} -
a^{jk}a^{j'k'}\ell_{x_j} v_{x_{j'}}v_{x_{k'}}
\Big)\nu^k d\Si \\
\ns\ds = \l \mathbb{E}\int_{\Si}
\sum_{j,k=1}^n\sum_{j',k'=1}^n\Big(
2a^{jk}a^{j'k'}\f_{x_{j'}}v_{x_j} v_{x_{k'}} -
a^{jk}a^{j'k'}\f_{x_j} v_{x_{j'}}v_{x_{k'}}
\Big)\nu^k d\Si
\\
\ns\ds = \l \mathbb{E}\int_{\Si}
\sum_{j,k=1}^n\sum_{j',k'=1}^n\Big(
2a^{jk}a^{j'k'}\f_{x_{j'}}\frac{\pa v}{\pa
\nu}\nu^j \frac{\pa v}{\pa \nu}\nu^{k'} -
a^{jk}a^{j'k'}\f_{x_j} \frac{\pa v}{\pa
\nu}\nu^{j'}\frac{\pa v}{\pa \nu}\nu^{k'}
\Big)\nu_k d\Si  \\
\ns\ds =\l \mathbb{E}\int_{\Si} \Big(
\sum_{j,k=1}^na^{jk}\nu^j \nu^k \Big)\Big(
\sum_{j',k'=1}^n a^{j'k'}\f_{x_{j'}}\nu^{k'}
\Big)\Big|\frac{\pa v}{\pa \nu}\Big|^2d\Si.
\end{array}
\end{equation}

For any $\tau\in (0,T_1)$ and $\tau'\in
(T_1',T)$, put
\begin{equation}\label{1.20-eq9}
Q_{\tau}^{\tau'}\=(\tau,\tau')\times G.
\end{equation}
Integrating (\ref{hyperbolic2}) in
$Q_{\tau}^{\tau'}$, taking expectation  and by
\eqref{coeffvt}, \eqref{bcoeffvtvi},
\eqref{8.21-eq1} and \eqref{8.21-eq5}, we obtain
that
\begin{eqnarray}\label{bhyperbolic3}
&&\ds \mathbb{E}\int_{Q_{\tau}^{\tau'}}
\theta\Big( -2\ell_t \hat v + 2\sum_{i=1}^n
a^{jk}\ell_{x_j}v_{x_k} + \Psi v \Big)
\Big[d\hat{\mathbf{z}} - \sum_{j,k=1}^n (a^{jk}\mathbf{z}_{x_j})_{x_k} dt\Big]dx\nonumber\\
&&\ds\q + \l \mathbb{E}\int_{\Si_0}\Big(
\sum_{j,k=1}^na^{jk}\nu^j \nu^k \Big)\Big(
\sum_{j',k'=1}^n a^{j'k'}\f_{x_{j'}}\nu^{k'}
\Big)\Big|\frac{\pa v}{\pa \nu}\Big|^2d\Si \nonumber\\
&&\ds \q + \mathbb{E}\int_{Q_{\tau}^{\tau'}}d
\[\ell_t\sum_{j,k=1}^n a^{jk}v_{x_j}
v_{x_k}+\ell_t\hat v^2-2\sum_{j,k=1}^na^{jk}\ell_{x_j} v_{x_k} \hat v - \Psi v\hat v+\(\cA\ell_t+\frac{\Psi_t}{2}\)v^2\]dx \\
&&\ds \geq c_0\l
\mathbb{E}\int_{Q_{\tau}^{\tau'}} \hat v^2 dxdt
+ \mathbb{E}\int_{Q_{\tau}^{\tau'}}
\sum_{j,k=1}^n c^{jk}v_{x_j}v_{x_k}dxdt +
\mathbb{E}\int_{Q_{\tau}^{\tau'}} \cB
v^2 dxdt \nonumber\\
&&\ds \q +
\mathbb{E}\!\int_{Q_{\tau}^{\tau'}}\!\!\[\ell_{t}\sum_{j,k=1}^n
\!a^{jk}(dv_{x_j})
(dv_{x_k})\!+\!\cA\ell_t(dv)^2\] dx\! +\!
\mathbb{E}\!\int_Q\!\Big(\!\!-2\ell_t\hat v\!
+\! 2\sum_{j,k=1}^n\!\!
a^{jk}\ell_{x_j}v_{x_k}\! +\! \Psi v\!\Big)^2
dxdt.\nonumber
\end{eqnarray}
Clearly,
\begin{eqnarray}\label{1.20-eq21}
&&\ds \mathbb{E}\int_{Q_{\tau}^{\tau'}}d
\[\ell_t\sum_{j,k=1}^n a^{jk}v_{x_j}
v_{x_k}+\ell_t\hat
v^2-2\sum_{j,k=1}^na^{jk}\ell_{x_j} v_{x_k} \hat
v - \Psi v\hat
v+\(\cA\ell_t+\frac{\Psi_t}{2}\)v^2\]dx\nonumber\\
&&\ds = \mathbb{E}\int_{G}
\[\ell_t\sum_{j,k=1}^n a^{jk}v_{x_j}(\tau')
v_{x_k}(\tau')+\ell_t\hat
v(\tau')^2-2\sum_{j,k=1}^na^{jk}\ell_{x_j}
v_{x_k}(\tau') \hat v(\tau') - \Psi v(\tau')\hat
v(\tau')\nonumber\\
&&\ds
\qq\qq+\(\cA\ell_t+\frac{\Psi_t}{2}\)v(\tau')^2\]dx
\\
&&\ds \q -\mathbb{E}\int_{G}
\[\ell_t\sum_{j,k=1}^n a^{jk}v_{x_j}(\tau)
v_{x_k}(\tau)+\ell_t\hat
v(\tau)^2-2\sum_{j,k=1}^na^{jk}\ell_{x_j}
v_{x_k}(\tau) \hat v(\tau) - \Psi v(\tau)\hat
v(\tau)\nonumber\\
&&\ds
\qq\qq+\(\cA\ell_t+\frac{\Psi_t}{2}\)v(\tau)^2\]dx\nonumber\\
&&\ds \leq C \l^3 \mathbb{E} \int_G \big\{\big[
\hat v(\tau)^2\! +\! |\nabla v(\tau)|^2 +\!
v(\tau)^2 \big] \!+ \!\big[  \hat v(\tau')^2 +\!
|\nabla v(\tau')|^2 +\!
v(\tau')^2\big]\big\}dx.\nonumber
\end{eqnarray}

From $\theta=e^{\ell}$ and \eqref{1.20-eq7}, we
know that there is a $\l_1>0$ such that for all
$\l>\l_1$,
\begin{equation}\label{1.20-eq11}
\l^3 \theta(\tau)\leq 1,\qq \l^3
\theta(\tau')\leq 1.
\end{equation}
Since $v=\theta \mathbf{z}$ and $\hat v=\theta
\hat{\mathbf{z}}$, it follows from
\eqref{1.20-eq11} that
\begin{equation}\label{1.20-eq10}
\begin{array}{ll}\ds
\l^3 \mathbb{E} \int_G \big\{\big[ \hat
v(\tau)^2 +\! |\nabla v(\tau)|^2 +\! v(\tau)^2
\big] \!+ \!\big[  \hat v(\tau')^2\!+\! |\nabla
v(\tau')|^2 +\!
v(\tau')^2\big]\big\}dx\\
\ns\ds \leq  \mathbb{E} \int_G \big\{\big[
\hat{\mathbf{z}}(\tau)^2 + |\nabla
\mathbf{z}(\tau)|^2 + \mathbf{z}(\tau)^2 \big] +
\big[  \hat{\mathbf{z}}(\tau')^2 + |\nabla
\mathbf{z}(\tau')|^2 +
\mathbf{z}(\tau')^2\big]\big\}dx.
\end{array}
\end{equation}
From \eqref{1.20-eq2}, \eqref{1.20-eq6} and
\eqref{1.20-eq9}, we obtain that $\L_2\subset
Q_\tau^{\tau'}$.  Thus,
\begin{equation}\label{1.2-eq1}
\begin{array}{ll}\ds
\l
\mathbb{E}\int_{Q_{\tau}^{\tau'}\setminus\L_2}
\hat v^2 dxdt =\l
\mathbb{E}\int_{Q_{\tau}^{\tau'}\setminus\L_2}
\big(\th\hat{\mathbf{z}} + \ell_t\th
\mathbf{z}\big)^2 dxdt\\
\ns\ds \leq 2\l
\mathbb{E}\int_{Q_{\tau}^{\tau'}\setminus\L_2}
\th^2\hat{\mathbf{z}}^2dxdt + 2\l^3
\mathbb{E}\int_{Q_{\tau}^{\tau'}\setminus\L_2}
(2t-T)^2\th^2 \mathbf{z}^2 dxdt
\end{array}
\end{equation}
and
\begin{equation}\label{1.2-eq2}
\begin{array}{ll}\ds
\mathbb{E}\int_{Q_{\tau}^{\tau'}\setminus\L_2}
\sum_{j,k=1}^n c^{jk}v_{x_j}v_{x_k}dxdt =
\mathbb{E}\int_{Q_{\tau}^{\tau'}\setminus\L_2}
\sum_{j,k=1}^n c^{jk}(\th \mathbf{z}_{x_j})(\th \mathbf{z}_{x_k})dxdt\\
\ns\ds =
\mathbb{E}\int_{Q_{\tau}^{\tau'}\setminus\L_2}
\sum_{j,k=1}^n c^{jk}\th^2(\l \f_{x_j}+
\mathbf{z}_{x_j})(\l \f_{x_k}+
\mathbf{z}_{x_k})dxdt\\
\ns\ds\leq C\l
\mathbb{E}\int_{Q_{\tau}^{\tau'}\setminus\L_2}\th^2
|\nabla \mathbf{z}|^2dxdt + C\l^3
\mathbb{E}\int_{Q_{\tau}^{\tau'}\setminus\L_2}\th^2
| \mathbf{z}|^2dxdt.
\end{array}
\end{equation}
Furthermore, it follows from \eqref{B1} that
\begin{equation}\label{1.2-eq3}
\begin{array}{ll}\ds
\mathbb{E}\int_{Q_{\tau}^{\tau'}\setminus\L_2}
\cB v^2 dxdt \leq C\l^3
\mathbb{E}\int_{Q_{\tau}^{\tau'}\setminus\L_2}
\th^2 \mathbf{z}^2 dxdt.
\end{array}
\end{equation}
Next, by \eqref{8.21-eq2} and \eqref{8.21-eq4},
we get that
\begin{equation}\label{8.21-eq6}
\begin{array}{ll}\ds
\mathbb{E}\!\int_{Q_{\tau}^{\tau'}\setminus\L_2}
\[\ell_{t}\sum_{j,k=1}^n  a^{jk}(dv_{x_j})
(dv_{x_k}) + \cA\ell_t(dv)^2\] dx\\
\ns\ds\leq C\l
\mathbb{E}\int_{Q_{\tau}^{\tau'}\setminus\L_2}\th^2
|\nabla \mathbf{z}|^2dxdt + C\l^3
\mathbb{E}\int_{Q_{\tau}^{\tau'}\setminus\L_2}\th^2
| \mathbf{z}|^2dxdt.
\end{array}
\end{equation}
From \eqref{1.20-eq2}, we know that $\th\leq
\exp(\l e^{R_0^2\mu/8})$ in
$Q_{\tau}^{\tau'}\setminus\L_2$. Consequently,
there exists $\l_2\geq\max\{\l_0,\l_1\}$  such
that for all $\l\geq\l_2$,
\begin{equation}\label{1.2-eq4}
C\l \max_{(x,t)\in
Q_{\tau}^{\tau'}\setminus\L_2} \th^2 \leq e^{\l
R_0^2 /3}, \qq C\l^3 \max_{(x,t)\in
Q_{\tau}^{\tau'}\setminus\L_2} \th^2 \leq e^{\l
R_0^2 /3}.
\end{equation}

It follows from \eqref{8.21-eq5} and
\eqref{1.2-eq1}--\eqref{1.2-eq4}  that
\begin{eqnarray}\label{1.20-eq11.1}
&& \l \mathbb{E}\int_{Q_{\tau}^{\tau'}} \hat
v^2 dxdt+\mathbb{E}\int_{Q_{\tau}^{\tau'}}
\sum_{j,k=1}^n c^{jk}v_{x_j}v_{x_k}dxdt +
\mathbb{E}\int_{Q_{\tau}^{\tau'}} \cB v^2 dxdt\nonumber\\
&&\q + \mathbb{E} \int_{Q_{\tau}^{\tau'}}
\[\ell_{t}\sum_{j,k=1}^n \!a^{jk}(dv_{x_j})
(dv_{x_k}) + \cA\ell_t(dv)^2\] dx\\
&& \geq \l \mathbb{E}\int_{\L_2}\! \hat v^2 dxdt
\!+ c_2\l \mathbb{E} \!\int_{\L_2}\! |\nabla
v|^2dxdt\!+ c_2\l^3 \mathbb{E} \!\int_{\L_2}\!
|v|^2dxdt \! - e^{\l R_0^2 /3}\mathbb{E}\!
\int_{Q}\! \big(|\hat{\mathbf{z}}|^2\!+ |\nabla
\mathbf{z}|^2\big)dxdt.\nonumber
\end{eqnarray}

Noting that $(\mathbf{z},\hat{\mathbf{z}})$
solves the equation \eqref{bsystem2}, we deduce
that
\begin{equation}\label{bhyperbolic4}
\begin{array}{ll}\ds
\mathbb{E}\int_{Q_{\tau}^{\tau'}} \theta \Big(
-2\ell_t \hat v + 2\sum_{i=1}^n
a^{jk}\ell_{x_j}v_{x_k} + \Psi v \Big)
\Big[d\hat{\mathbf{z}} - \sum_{j,k=1}^n (a^{jk}\mathbf{z}_{x_j})_{x_k} dt\Big]dx \\
\ns\ds= \mathbb{E}\int_{Q_{\tau}^{\tau'}}
\theta\Big(\! -2\ell_t \hat v + 2\sum_{i=1}^n
a^{jk}\ell_{x_j}v_{x_k}\! + \Psi v \Big)\big[-a_1\!\cd\!\nabla \mathbf{z} + \big(-\div a_1\!+\!a_2\!-\!a_3a_5\big) \mathbf{z} \big]dxdt\\
\ns\ds\leq \mathbb{E}\int_{Q_{\tau}^{\tau'}}
\theta\Big(\! -2\ell_t \hat v + 2\sum_{i=1}^n
a^{jk}\ell_{x_j}v_{x_k}\! + \Psi v \Big)^2dxdt +
r_2\mathbb{E}\int_{Q_{\tau}^{\tau'}}(|\n\mathbf{z}|^2
+ \mathbf{z}^2)dxdt.
\end{array}
\end{equation}
Combing \eqref{bhyperbolic3}, \eqref{1.20-eq10},
\eqref{1.20-eq11.1} and \eqref{bhyperbolic4}, we
conclude that there is a $\l_3
\geq\max\{\l_2,Cr_5+1\}$ such that for any $\l
\geq \l_3$, one has that
\begin{equation}\label{1.20-eq11.2}
\begin{array}{ll}\ds
\mathbb{E} \int_{\L_1}
\th^2(|\hat{\mathbf{z}}|^2+|\nabla
\mathbf{z}|^2)dxdt+  \mathbb{E}
\int_{\L_1} \th^2 |\mathbf{z}|^2dxdt\\
\ns\ds \leq \!  C \Big[e^{\l R_0^2 /3}\mathbb{E}
\int_{Q}\big(|\hat{\mathbf{z}}|^2+|\nabla
\mathbf{z}|^2\big)dxdt + e^{\l R_0^2
/3}\mathbb{E} \int_{Q} |\mathbf{z}|^2dxdt  +
e^{\l R_1^2} \mathbb{E} \int_{\Si_0}
\Big|\frac{\pa \mathbf{z}}{\pa
\nu}\Big|^2d\Si \\
\ns\ds \qq  + \mathbb{E} \int_G \big(
\hat{\mathbf{z}}(\tau)^2 + |\nabla
\mathbf{z}(\tau)|^2 + \mathbf{z}(\tau)^2 +
\hat{\mathbf{z}}(\tau')^2 + |\nabla
\mathbf{z}(\tau')|^2 +
\mathbf{z}(\tau')^2\big)dx\Big].
\end{array}
\end{equation}
Integrating \eqref{1.20-eq11.2} with respect to
$\tau$ and $\tau'$ on $[T_2,T_1]$ and
$[T_1',T_2']$, respectively, we get that
\begin{equation}\label{1.20-eq12}
\begin{array}{ll}\ds
\mathbb{E} \int_{\L_1}
\th^2(|\hat{\mathbf{z}}|^2+|\nabla
\mathbf{z}|^2)dxdt+  \mathbb{E}
\int_{\L_1} \th^2 |\mathbf{z}|^2dxdt\\
\ns\ds \leq   C e^{\l R_0^2 /3}\mathbb{E}
\int_{Q} \big(|\hat{\mathbf{z}}|^2+|\nabla
\mathbf{z}|^2\big)dxdt + C e^{\l R_0^2
/3}\mathbb{E} \int_{Q}
|\mathbf{z}|^2dxdt  \\
\ns\ds \q +  Ce^{\l R_1^2}\mathbb{E}
\int_{\Si_0} \Big|\frac{\pa \mathbf{z}}{\pa
\nu}\Big|^2d\Si + C\mathbb{E} \int_Q \big(
\hat{\mathbf{z}}^2 + |\nabla \mathbf{z}|^2 +
\mathbf{z}^2 \big)dxdt.
\end{array}
\end{equation}
From \eqref{1.20-eq2}, we obtain that
\begin{equation}\label{1.20-eq13}
\begin{array}{ll}\ds
\mathbb{E} \int_{\L_1}\th^2
(|\hat{\mathbf{z}}|^2+|\nabla \mathbf{z}|^2
+|\mathbf{z}|^2)dxdt \3n&\ds\geq \mathbb{E}
\int_{\L_0} \th^2(|\hat{\mathbf{z}}|^2+|\nabla
\mathbf{z}|^2+|\mathbf{z}|^2)dxdt\\
\ns&\ds \geq e^{\l R_0^2 /2}\mathbb{E}
\int_{Q_0} (|\hat{\mathbf{z}}|^2+|\nabla
\mathbf{z}|^2+|\mathbf{z}|^2)dxdt.
\end{array}
\end{equation}
Combing \eqref{1.20-eq12} and \eqref{1.20-eq13},
we arrive at
\begin{equation}\label{1.20-eq14}
\begin{array}{ll}\ds
 \mathbb{E} \int_{Q_0}
(|\hat{\mathbf{z}}|^2+|\nabla
\mathbf{z}|^2+|\mathbf{z}|^2)dxdt\\
\ns\ds \leq   C \[e^{-\l R_0^2 /6}\mathbb{E}
\int_{Q} (|\hat{\mathbf{z}}|^2+|\nabla
\mathbf{z}|^2+|\mathbf{z}|^2)dxdt + \l^2e^{-\l
R_0^2 /6}\mathbb{E} \int_{Q}
|\mathbf{z}|^2dxdt \\
\ns\ds \qq + e^{\l R_1^2}\mathbb{E} \int_{\Si_0}
\Big|\frac{\pa \mathbf{z}}{\pa \nu}\Big|^2d\Si +
\mathbb{E} e^{-\l R_0^2 /2}\int_Q \big(
\hat{\mathbf{z}}^2 + |\nabla \mathbf{z}|^2 +
\mathbf{z}^2 \big)dxdt\].
\end{array}
\end{equation}
By standard energy estimate of the equation
\eqref{bsystem2}, we have that
\begin{equation}\label{1.20-eq15}
\begin{array}{ll}\ds
|(\mathbf{z}_0,\hat{\mathbf{z}}_0)|_{H_0^1(G)\times
L^2(G)}^2 \leq Ce^{Cr_4} \mathbb{E} \int_{Q_0}
\big(|\hat{\mathbf{z}}|^2+|\nabla
\mathbf{z}|^2+|\mathbf{z}|^2\big)dxdt
\end{array}
\end{equation}
and that
\begin{equation}\label{1.20-eq16}
\begin{array}{ll}\ds
|(\mathbf{z}_0,\hat{\mathbf{z}}_0)|_{H_0^1(G)\times
L^2(G)}^2 \geq Ce^{-Cr_4} \mathbb{E} \int_{Q}
\big(|\hat{\mathbf{z}}|^2+|\nabla
\mathbf{z}|^2+|\mathbf{z}|^2\big)dxdt.
\end{array}
\end{equation}

It follows from
\eqref{1.20-eq14}--\eqref{1.20-eq16} that
\begin{equation}\label{1.20-eq14-1}
\begin{array}{ll}\ds
|(\mathbf{z}_0,\hat{\mathbf{z}}_0)|_{H_0^1(G)\times
L^2(G)}^2 \leq   C e^{Cr_4}e^{-\l R_0^2
/6}|(\mathbf{z}_0,\hat{\mathbf{z}}_0)|_{H_0^1(G)\times
L^2(G)}^2   +  C e^{\l R_1^2}\mathbb{E}
\int_{\Si_0} \Big|\frac{\pa \mathbf{z}}{\pa
\nu}\Big|^2d\Si.
\end{array}
\end{equation}
Let us choose $\l_4\geq \l_3$ such that $C
e^{Cr_4}e^{-\l_4 R_0^2 /6}<1$. Then, for all
$\l\geq \l_4$, we have that
\begin{equation}\label{bhyperbolic8}
\begin{array}{ll}
\ds\q
|(\mathbf{z}_0,\hat{\mathbf{z}}_0)|_{H_0^1(G)\times
L^2(G)}^2 \leq Ce^{\l R_1^2}
\mathbb{E}\int_{\Si_0} \Big|\frac{\pa
\mathbf{z}}{\pa \nu}\Big|^2d\Si.
\end{array}
\end{equation}
This leads to the inequality \eqref{exact ob
est} immediately.
\endpf

\begin{remark}
It follows from Proposition \ref{prop-hid1}
that
$$
\Big|\frac{\pa
\mathbf{z}}{\pa\nu}\Big|_{L^2_\dbF(0,T;L^2(\G_0))}
\leq
C\big(|\mathbf{z}_0|_{H^{1}_0(G)}+|\hat{\mathbf{z}}_0|_{L^2(G)}\big).
$$
This, together with Theorem \ref{th-Ob}, implies
that
$$
\frac{1}{C}\big(|\mathbf{z}_0|_{H^{1}_0(G)}+|\hat{\mathbf{z}}_0|_{L^2(G)}\big)\leq
\Big|\frac{\pa
\mathbf{z}}{\pa\nu}\Big|_{L^2_\dbF(0,T;L^2(\G_0))}
\leq
C\big(|\mathbf{z}_0|_{H^{1}_0(G)}+|\hat{\mathbf{z}}_0|_{L^2(G)}\big).
$$
Therefore, we can defined a new norm
$|\!|\cd|\!|$ on $H^{1}_0(G)\times L^2(G)$,
which is equivalent to the norm
$|\cd|_{H^{1}_0(G)\times L^2(G)}$ as
follows:\vspace{-0.1cm}
$$
\begin{array}{ll}\ds
|\!|(\xi,\eta)|\!| = \Big|\frac{\pa
\mathbf{z}}{\pa\nu}\Big|_{L^2_\dbF(0,T;L^2(\G_0))}
\q \mbox{ for }(\xi,\eta)\in H^{1}_0(G)\times
L^2(G),
\end{array}
$$
where $(\mathbf{z},\hat{\mathbf{z}})$ is the
solution to \eqref{bsystem2} and
$(\mathbf{z}_0,\hat{\mathbf{z}}_0)=(\xi,\eta)$.
This is the starting point of the duality
argument to the proof of the controllability
result via the observability estimate (See
\cite{Lions1} for the details for wave
equations).
\end{remark}
%


\section{Proofs of the main results}
\label{sec-con}


This section is addressed to proving our main
results in this paper, i.e., Theorems
\ref{th-non-con}--\ref{th-non-con1}.

{\it Proof of  Theorem \ref{exact th}}\,: It
follows from Propositions \ref{prop-exact} and
\ref{prop-dual1}, and Theorem \ref{th-Ob}
immediately.
\endpf

Before proving Theorems \ref{th-non-con} and
\ref{th-non-con1}, we recall the following known
result (\cite[Lemma 2.1]{Peng}).

\begin{lemma}\label{lm1}
There is a random variable $\xi\in
L^2_{\cF_T}(\Om)$ such that it is impossible to
find $(\varrho_1,\varrho_2)\in
L^2_\dbF(0,T)\times C_\dbF([0,T];L^2(\Om))$ and
$\a\in \dbR$ satisfying
$$
\xi = \a + \int_0^T \varrho_1(t)dt + \int_0^T
\varrho_2(t)dW(t).
$$
\end{lemma}

\ms

We are now in a position to prove Theorems
\ref{th-non-con} and \ref{th-non-con1}.

\ms

{\it Proof of Theorem \ref{th-non-con}}\,: We
use the contradiction argument. Choose $\psi\in
H_0^1(G)$ satisfying $|\psi|_{L^2(G)}=1$ and let
$\tilde y_0=\xi \psi$, where $\xi$ is given in
Lemma \ref{lm1}. Assume that \eqref{system1} was
exactly controllable. Then, for any $y_0\in
L^2(G)$, we would find a triple of controls
$(g_1,g_2,h)\in L^2_{\dbF}(0,T;H^{-1}(G))\times
L^2_{\dbF}(0,T;H^{-1}(G))\times
L^2_\dbF(0,T;L^2(\G_{0}))$ such that the
corresponding solution $ y\in
C_{\dbF}([0,T];L^2(\Om;L^2(G)))\cap
C_{\dbF}^1([0,T];L^2(\Om;H^{-1}(G))) $ to the
equation \eqref{system1} satisfies that
$y(T)=\tilde y_0$. Clearly,
$$
\int_G \tilde y_0 \psi dx-\int_G  y_0 \psi dx =
\int_0^T\langle
y_t,\psi\rangle_{H^{-1}(G),H_0^1(G)}dt,
$$
which leads to
$$
\xi = \int_G  y_0 \psi dx + \int_0^T\langle
y_t,\psi\rangle_{H^{-1}(G),H_0^1(G)}dt.
$$
This contradicts Lemma \ref{lm1}.
\endpf

\ms

{\it Proof of Theorem \ref{th-non-con1}}\,: Let
us employ the contradiction argument, and divide
the proof into three cases.

\ms

\textbf{Case 1) $a_4\in
C_{\dbF}([0,T];L^\infty(\Omega))$ and $f$ is
supported in $G_0$}. Since $G_0\subset G$ is an
open subset and $G\setminus
\overline{G_0}\neq\emptyset$, we can find a
$\rho\in C_0^\infty(G\setminus G_0)$ satisfying
$|\rho|_{L^2(G)}=1$.

Assume that \eqref{system2} was exactly
controllable. Then, for $(y_0,\hat y_0)=(0,0)$,
one could find controls $(f,g,h)\in
L_{\dbF}^2(0,T;L^2(G))\times
L_{\dbF}^2(0,T;H^{-1}(G))\times L_{\dbF}^2(0,T;
L^2(\G_{0}))$ with $\supp f\subset G_0$, a.e.
$(t,\omega)\in (0,T)\times\Omega$ such that the
corresponding solution to \eqref{system2}
fulfills $(y(T),\hat y(T)) =(\rho\xi,0)$, where
$\xi$ is given in Lemma \ref{lm1}.
Thus,\vspace{-0.1cm}
\begin{equation}\label{8.27-eq2}
\rho\xi = \int_0^T \hat y dt + \int_0^T
(a_4y+f)dW(t).
\end{equation}
Multiplying both sides of \eqref{8.27-eq2} by
$\rho$ and integrating it in $G$, we get that
\begin{equation}\label{8.27-eq3}
\xi = \int_0^T \langle\hat
y,\rho\rangle_{H^{-1}(G),H_0^1(G)}dt + \int_0^T
\langle a_4y,\rho\rangle_{L^2(G)} dW(t).
\end{equation}
Since the pair $(y,\hat y)\in
C_\dbF([0,T];L^2(\Om;L^2(G)))\times
C_\dbF([0,T];L^2(\Om;H^{-1}(G)))$ solves
\eqref{system2}, then
$ \langle\hat y,\rho\rangle_{H^{-1}(G),H_0^1(G)}
\in C_\dbF([0,T];L^2(\Om))$ and $ \langle
y,\rho\rangle_{L^2(G)}$ $ \in
C_\dbF([0,T];L^2(\Om)), $
which, together with \eqref{8.27-eq3},
contradicts Lemma \ref{lm1}.

\ms

\textbf{Case 2) $a_3\in
C_{\dbF}([0,T];L^\infty(\Omega))$ and $g$ is
supported in $G_0$}. Choose $\rho$ as in Case
1).

If \eqref{system2} was exactly controllable,
then, for $(y_0,\hat y_0)=(0,0)$, one can find
controls $(f,g,h)\in
L_{\dbF}^2(0,T;L^2(G))\times
L_{\dbF}^2(0,T;H^{-1}(G))\times L_{\dbF}^2(0,T;
L^2(\G_{0}))$ with $\supp g\subset G_0$, a.e.
$(t,\omega)\in (0,T)\times\Omega$  such that the
corresponding solution of \eqref{system2}
fulfills $(y(T),\hat y(T)) =(0,\xi)$.

It is clear that $(\phi,\hat\psi)=(\rho y, \rho
\hat y)$ solves the following  equation:
\begin{equation}\label{8.27-eq1}
\left\{ \begin{array}{ll}
\ds d\phi= \hat \phi dt+(a_4\phi+\rho f)dW(t) &\mbox{ in }Q,\\
\ns\ds d\hat \phi-
\sum_{j,k=1}^n(a^{jk}\phi_{x_j})_{x_k}dt= \zeta
dt + a_3\phi
dW(t)&\mbox{ in }Q,\\
\ns\ds \phi= 0 &\mbox{ on }\Si ,\\
\ns\ds \hat \phi=0&\mbox{ on }\Si,\\
\ns\ds \phi(0)=0,\q \hat \phi(0)=0&\mbox{ in }G,
\end{array}
\right.
\end{equation}
where $\ds\zeta =
\sum_{j,k=1}^n[(a^{jk}\rho_{x_j} y)_{x_k} +
a^{jk}y_{x_j}\rho_{x_k}]+\rho a_1\cdot\nabla
y+\rho a_2y$. Further, we have $\phi(T)=0$ and
$\hat\phi(T)=\rho\xi$. Noting that
$(\phi,\hat\phi)$ is the weak solution to
\eqref{8.27-eq1}, we see that
$$
\begin{array}{ll}\ds
\big\langle
\rho\xi,\rho\big\rangle_{H^{-2}(G),H_0^2(G)}\\
\ns\ds = \int_0^T\!
\[\Big\langle\sum_{j,k=1}^n\!(a^{jk}\phi_{x_j})_{x_k},\rho\Big\rangle_{H^{-2}(G),H_0^2(G)}\! + \big\langle\zeta,\rho\big\rangle_{H^{-1}(G),H_0^1(G)}
\]dt +\! \int_0^T\!\big\langle
a_3\phi,\rho\big\rangle_{L^2(G)}dW(t),
\end{array}
$$
which implies that
\begin{equation}\label{8.27-eq4}
\begin{array}{ll}\ds
\xi \!\3n&\ds= \!\int_0^T\!
\[\Big\langle\!\sum_{j,k=1}^n\!(a^{jk}\phi_{x_j})_{x_k},\rho\Big\rangle_{H^{-2}(G),H_0^2(G)} \!+ \big\langle\zeta,\rho\big\rangle_{H^{-1}(G),H_0^1(G)}
\]dt\! +\! \int_0^T\!\big\langle
a_3\phi,\rho\big\rangle_{L^2(G)}dW(t).
\end{array}
\end{equation}
Since $(\phi,\hat \phi)\in
C_\dbF([0,T];L^2(\Om;L^2(G)))\times
C_\dbF([0,T];L^2(\Om;H^{-1}(G)))$, then
$$
\Big\langle\sum_{j,k=1}^n(a^{jk}\phi_{x_j})_{x_k},\rho\Big\rangle_{H^{-2}(G),H_0^2(G)}
+
\big\langle\zeta,\rho\big\rangle_{H^{-1}(G),H_0^1(G)}
\in L^2_\dbF(0,T)
$$
and\vspace{-0.21cm}
$$
\langle a_3\phi,\rho\rangle_{L^2(G)} \in
C_\dbF([0,T];L^2(\Om)).
$$
These, together with \eqref{8.27-eq4},
contradict Lemma \ref{lm1}.

\ms

\textbf{Case 3) $h=0$}. Assume that the system
\eqref{system2} was exactly controllable.
Similar to the proof of Proposition
\ref{prop-dual1}, we could deduce that, for any
$(z^T,\hat z^T)\in
L^2_{\cF_T}(\Om;H_0^1(G))\times
L^2_{\cF_T}(\Om;L^2(G))$, the solution
$(z,Z,\hat z,\widehat Z)$ to \eqref{bsystem1}
(with $\tau=0$ and $(z(T),\hat z(T))=(z^T,\hat
z^T)$) satisfies
\begin{equation}\label{9.3-eq6}
|(z^T,\hat
z^T)|_{L^2_{\cF_T}(\Om;H^{1}_0(G))\times
L^2_{\cF_T}(\Om;L^2(G))} \leq C \(
|Z|_{L^2_\dbF(0,T;H^1_0(G))}+|\widehat
Z|_{L^2_\dbF(0,T;L^2(G))}\).
\end{equation}
For any nonzero $(\eta_0, \eta_1)\in
H_0^1(G)\times L^2(G)$, let us consider the
following random wave equation:
\begin{equation}\label{9.3-eq7}
\left\{
\begin{array}{ll}
\ds \eta_{tt} - \sum_{j,k=1}^n(a^{jk}\eta_{x_j})_{x_k} dt = b_1\cd\nabla \eta + b_2 \eta &\mbox{ in } (0,T)\times G,\\
\ns\ds z = 0 &\mbox{ on } (0,\tau)\times \G,\\
\ns\ds \eta(0) =\eta_0,\q \eta_t(0) = \eta_1
&\mbox{ in } G.
\end{array}
\right.
\end{equation}
Clearly, $(\eta,0,\eta_t,0)$ solves
\eqref{bsystem1} with the final datum $(z^T,\hat
z^T)=(\eta(T),\hat\eta(T))$, a contradiction to
the inequality \eqref{9.3-eq6}.
\endpf


\section{Further comments and open
problems}\label{sec-com}

In this paper, we obtain the exact controllability
of the system \eqref{system2} with one boundary
control and two internal controls.  It is
natural to consider the exact controllability
problem for stochastic wave-like equations with
three internal controls:
\begin{equation*}\label{system6}
\left\{
\begin{array}{ll}
\ds dy= \hat y dt+(a_4y+f)dW(t) &\mbox{ in }Q,\\
\ns\ds d\hat
y-\sum_{j,k=1}^n(a^{jk}y_{x_j})_{x_k}dt=(a_1 \cd
\nabla y + a_2 y + a_5 g+\chi_{G_0}h)dt +
(a_3y+g)
dW(t)&\mbox{ in }Q,\\
\ns\ds y= 0 &\mbox{ on }\Si ,\\
\ns\ds y(0)=y_0,\q \hat y(0)=\hat y_0&\mbox{ in
}G.
\end{array}
\right.
\end{equation*}
Here \vspace{-0.1cm}
$$
G_0\=\{x\in G\,|\, \dist(x,\G_0)\leq
\d\}\vspace{-0.1cm}
$$
for a $\d>0$, $(y_0,\hat y_0)\in H_0^1(G)\times
L^2(G)$, $f \in L^2_\dbF(0,T;H_0^1(G))$, $g\in
L^2_\dbF(0,T;L^2(G))$ and $h\in
L^2_\dbF(0,T;L^2(G_0))$ are the controls.

By a duality argument, one only need to show the
following observability estimate:
\begin{equation}\label{4.1-eq1}
\begin{array}{ll}\ds
|z^T,\hat z^T|_{L^2_{\cF_T}(\Om;L^2(G))\times
L^2_{\cF_T}(\Om;H^{-1}(G))}\\
\ns\ds\leq
C\(|\chi_{G_0}z|_{L^2_\dbF(0,T;L^2(\G_0))}+|Z|_{L^2_\dbF(0,T;L^2(G))}+|\widehat
Z|_{L^2_\dbF(0,T;H^{-1}(G))}\),
\end{array}
\end{equation}
where $(z,Z,\hat z,\widehat Z)$ solves
\eqref{bsystem1} with $\tau=T$ and the final
datum $(z^T,\hat z^T)$. Following  \cite{FLLZ},
we can prove that \eqref{4.1-eq1} holds. Details
are too lengthy to be presented here.

\vspace{0.2cm}

There are many open problems related to the topic
of this paper. We shall list below some of them which, in our opinion, are particularly
interesting:

\begin{itemize}

\item {\bf Null controllability for
stochastic wave-like equations}

In this paper,  the exact controllability for
stochastic wave-like equations are presented. As
immediate consequences, we can obtain the null
and approximate controllability for the same
system. However, in order to show these two
results, there seems no reason to use three
controls. By Theorem \ref{th-non-con1}, it is
shown that only one control applied in the
diffusion term is not enough. However, inspired
by the result in \cite{Lu0}, we believe that one
boundary control is enough for the null and
approximate controllability of \eqref{system1}
and \eqref{system2}. Unfortunately, some
essential difficulties appear when we try to
prove it, following  the method in the present
paper. For example,  for the null
controllability, we should prove the following
inequality for the solution to \eqref{bsystem1}:
$$
|z(0)|^2_{H_0^1(G)} + |\hat z(0)|^2_{L^2(G)}
\leq C\int_0^T\int_{\G_0} \Big| \frac{\pa
z}{\pa\nu} \Big|^2 d\G dt.
$$
However, if we utilize the method in this paper,
we only get
$$
\begin{array}{ll}\ds
|z(0)|^2_{H_0^1(G)} +  |\hat
z(0)|^2_{L^2(G)}\leq C\(\int_0^T\int_{\G_0}
\Big| \frac{\pa z}{\pa\nu} \Big|^2 d\G dt +
\int_0^T |Z|_{H_0^1(G)}^2 + \int_0^T |\widehat
Z|_{L^2(G)}^2 dt\).
\end{array}
$$
There are two additional terms containing $Z$
and $\widehat Z$ in the right hand side of the above inequality. These
terms come from the fact that, in the Carleman
estimate, we regard $Z$ and $\widehat Z$ simply as
nonhomogeneous terms rather than  part of the
solution. Therefore, we
believe that one should introduce some new
technique, for example, a Carleman estimate in
which the fact that $Z$ and $\widehat Z$ are  part of
the solution is essentially used, to get rid of the
additional terms containing $Z$ and $\widehat
Z$. However, we do not know how to achieve this
goal at this moment.

\item{\bf Exact controllability for stochastic wave-like equations
with less restrictive condition}

In this paper, we get the exact controllability
of the system \eqref{system2} for $\G_0$ given
by \eqref{def gamma0}. It is well known that a
sharp sufficient condition for exact
controllability of deterministic wave equations
with time invariant coefficients is that the
triple $(G, \G_0, T)$ satisfies the Geometric
Control Condition introduced in
\cite{Bardos-Lebeau-Rauch1}. It would be quite
interesting and challenging to extend this
result to the stochastic setting, but it seems
that there are lots of things to be done before
solving this problem. Indeed,  the propagation
of singularities for stochastic partial
differential equations, at least, for stochastic
hyperbolic equations, should be established.
However, as far as we know, this topic is
completely open.

\item{\bf Exact controllability for stochastic wave-like equations
with more regular controls}

In this paper, we get the exact controllability
of the system \eqref{system2} a triple
$(f,g,h)$, where $g\in L^2_\dbF(0,T;H^{-1}(G))$,
which is very irregular.  It is very interesting
to see whether \eqref{system2} is exactly
controllable when $g\in L^2_\dbF(0,T;L^2(G))$.
By duality argument, one can show that this is
equivalent to the following observability
estimate:
\begin{equation}\label{4.1-eq1.1}
\begin{array}{ll}\ds
|(z^T,\hat
z^T)|_{L^2_{\cF_T}(\Om;H_0^1(G))\times
L^2_{\cF_T}(\Om;L^2(G))}\\
\ns\ds\leq C\(\Big|\frac{\pa
z}{\pa\nu}\Big|_{L^2_\dbF(0,T;L^2(\G_0))}+|Z|_{L^2_\dbF(0,T;L^2(G))}+|\widehat
Z|_{L^2_\dbF(0,T;L^2(G))}\),
\end{array}
\end{equation}
where $(z,Z,\hat z,\widehat Z)$ is the solution
to \eqref{bsystem1} with $\tau=T$ and final
datum $(z^T,\hat z^T)$. By Lemma
\ref{hyperbolic1}, we can prove that the
inequality \eqref{4.1-eq1.1} holds if the term
$|Z|_{L^2_\dbF(0,T;L^2(G))}$ is replaced by
$|Z|_{L^2_\dbF(0,T;H^1_0(G))}$. However, we do
not know whether \eqref{4.1-eq1.1} is true or
not.

\end{itemize}



\end{document}